\documentclass[preprint,11pt]{amsart}
\usepackage[margin=1in]{geometry}              % to include figures
\usepackage{amsmath}               % great math stuff
\usepackage{amsfonts}              % for blackboard bold, etc
\usepackage{amsthm}
\usepackage{enumerate}

\numberwithin{equation}{section}

\newtheorem{thm}{Theorem}[section]
\newtheorem{defn}[thm]{Definition}
\newtheorem{lem}[thm]{Lemma}
\newtheorem{prop}[thm]{Proposition}
\newtheorem{cor}[thm]{Corollary}

\newtheorem{Rem}[thm]{Remark}

\def\XXint#1#2#3{{\setbox0=\hbox{$#1{#2#3}{\int}$}
     \vcenter{\hbox{$#2#3$}}\kern-.5\wd0}}

\begin{document}

\nocite{*}

\title[Regularity for obstacle problems]{Regularity for obstacle problems to anisotropic parabolic equations }


\author{Hamid EL Bahja}

\address{Hamid EL Bahja, The African Institute for Mathematical Sciences, Research and Innovation Centre, Rwanda}
\email{hamidsm88@gmail.com, helbahja@aims.ac.za}


\begin{abstract}
Following Dibenedetto's intrinsic scaling method, we prove local H\"older continuity of weak solutions to obstacle problems related to some anisotropic parabolic equations under the condition for which only H\"older's continuity of the obstacle is known.
\end{abstract}

\maketitle

\section{Introduction}
In this work, we consider the regularity issue for a class of anisotropic parabolic equations of the form
\begin{equation}
    u_{t}-\sum_{i=1}^{N}\frac{\partial}{\partial x_{i}}\left( \left|\frac{\partial u}{\partial x_{i}} \right|^{p_{i}-2}\frac{\partial u}{\partial x_{i}}\right)=0~~\text{     in}~\Omega_{T},
    \label{eq:main}
\end{equation}
where $\Omega_{T}\equiv \Omega\times(0,T]$, $\Omega$ is a bounded domain in $\mathbb{R}^{N}$, $N\geq2$, $T>0$, with the exponents $p_{i}\geq2$ for all $i=1,..,N$. The solutions to (1.1) are subject to an obstacle constraint of the form $u\geq \phi$ in $\Omega_{T}$, with $\phi$ being H\"older continuous. In recent decades, there has been growing interest in these types of equations because of their interesting feature of anisotropic diffusion with orthotropic structure where the diffusion rates diﬀer according to the direction $x_i$. Besides its inherent mathematical interest, they emerge for instance, from the mathematical description of the dynamics of fluids with different conductivities in different directions. This is important for modeling diffusion in materials that have a specific structure, such as wood, bone, or composite materials. For example, in bone, the diffusion of minerals is faster along the long axis of the bone than across the short axis. This is because the bone is made up of a network of interconnected canals, and the canals are oriented along the long axis of the bone. For more examples, see \cite{pie,mur,ant} and references therein.

In order to state our main result, we need to briefly describe the by-now classical approach to the regularity of solutions to the degenerate parabolic $p$-Laplace operator as first introduced by Dibenedetto \cite{dib1,dib2}. The latter realized that the poor structure of PDEs with quasi-linear parabolic operators should be taken into account. By including the singularity/degeneracy of the equation in a suitable geometry, we can derive integral inequalities for the "right" cylinders, which suggests that the PDE behaves in a specific way in its own geometry, and then the continuity of the solution at a point follows from showing that the oscillation converges to zero as a sequence of nested cylinders shrinks to the point.

Equation (1.1) with $p_{i}=p$ is reduced to the standard parabolic $p-$Laplacian type of equations whose regularity properties are well studied \cite{dib1,lad,ser,tru}. However, the theory for the obstacle case is not yet complete. H\"older continuity for a class of parabolic quasi-linear obstacle problems is presented in \cite{bog1,kor, kuu} and references therein, and for obstacle problems to porous medium type equations has been treated in the recent papers \cite{bog2,mor,cho}.  

When $p_{i}'s$ are potentially different, the regularity of local solutions to (1.1) has been studied by several authors. We refer, e.g., to  \cite{dibgian,cian,duz} for results on local continuity of solutions to (1.1), to \cite{elb24} for the local boundedness of the corresponding nonlocal equation, and to \cite{bou,rak} for higher regularity properties. Nevertheless, despite the previously mentioned results, the regularity theory for obstacle problems related to elliptic/parabolic anisotropic equations is largely unknown. The latter is precisely the aim of this paper, where we will prove the following local continuity result.

\begin{thm}
    Under the assumptions
    \[
    2 < p_i < \bar p \left( 1 + \frac{2}{N} \right), \qquad
    \bar p := \left( \frac{1}{N} \sum_{i=1}^N \frac{1}{p_i} \right)^{-1} < N,
    \qquad
    \bar p - p^- \le \frac{\bar p (p^- - 2)}{N},
    \]
    where \(p^- := \min_i p_i\), and the obstacle \(\phi \in C^{0;\beta,\beta/2}(\Omega_T)\) for some \(\beta \in (0,1)\), any local weak solution to the obstacle problem related to (1.1) in the sense of Definition 2.3 is locally H\"older continuous.
\end{thm}

\begin{Rem}
    The closeness condition \(\bar p - p^- \le \frac{\bar p (p^- - 2)}{N}\) ensures the De Giorgi iteration in Lemma 5.3 closes without a scale-dependent exponential penalty. This condition is similar to the one used in \cite{dibgian}. In particular, if \(\bar p - p^- \ll 1\), the condition is automatically satisfied, since the left-hand side becomes arbitrarily small while the right-hand side remains strictly positive (as \(p^- > 2\)).
\end{Rem}

To prove our result, we use a similar strategy to that used in \cite{bog2,mor} for nonnegative obstacles, which relies on energy estimations for truncations of the solution and on a De Giorgi-type iteration argument. The basic idea is to construct a sequence of cylinders shrinking to a common vertex. Within each of these cylinders, we consider two measure-theoretic alternatives, which we will call the first and second alternatives. In both alternatives, the solution is bounded away from one of its extrema in a quantifiable way pointwise almost everywhere in a smaller cylinder. The derivation of the energy estimates for the obstacle problem is more complicated than for the obstacle-free case because the solution is not differentiable in time. This prevents us from using the solution itself as a comparison map. To overcome this difficulty, we use a mollification argument in time to exploit the weak formulation of the obstacle problem in the sense of Definition 2.3. However, since our equation exhibits degenerate anisotropic scaling behavior, we will work in cylinders that respect the intrinsic geometry of the equation. 

This paper is organized as follows. In Section 2, we give the definition of weak solutions related to our obstacle problem and introduce some fundamental analytic tools. In Section 3, we establish local energy and logarithmic estimates. In Section 4, we use the energy estimates to prove the local boundedness of the weak solutions. Finally, in Section 5, we prove Theorem 1.1.

\section{Definitions and technical tools}
In what follows, we recall some definitions and basic properties of the anisotropic Lebesgue-Sobolev spaces. Then, for exponents $\{p_{i}\}_{i=1}^{N}\geq 1$ we introduce the anisotropic space
\begin{equation*}
    W^{1,p_{i}}(\Omega):=\{u\in L^{p_{i}}(\Omega),~\frac{\partial u}{\partial x_{i}}\in L^{p_{i}}(\Omega)\},
\end{equation*}
which is Banach space under the norm 
\begin{equation*}
    \|u\|_{W^{1,p_{i}}(\Omega)}:=\|u\|_{L^{p_{i}(\Omega)}}+\left\|\frac{\partial u}{\partial x_{i}}\right\|_{L^{p_{i}}(\Omega)}.
\end{equation*}
Also, $W_{0}^{1,p_i}(\Omega)$ denotes the closure of $C^{1}_{0}(\Omega)$ in $W^{1,p_i}(\Omega)$ for all $i=1,..,N$. Next, for a multi-index $\Vec{p}=(p_{1},..,p_{N})$ we put
\begin{equation*}
W^{1,\Vec{p}}(\Omega)= \bigcap_{i=1}^{N}W^{1,p_i}(\Omega),~\text{and}~L^{\Vec{p}}(0,T; W^{1,\Vec{p}}(\Omega))=\bigcap_{i=1}^{N}L^{p_i}(0,T;W^{1,p_i}(\Omega)),
\end{equation*}
with
\begin{equation*}
    \|u\|_{W^{1,\Vec{p}}(\Omega)}=\sum_{i=1}^{N}\|u\|_{W^{1,p_i}(\Omega)},~\text{and}~\|u\|_{L^{\Vec{p}}(0,T; W^{1,\Vec{p}}(\Omega))}=\sum_{i=1}^{N}\|u\|_{L^{p_i}(0,T;W^{1,p_i}(\Omega))}.
\end{equation*}
We state now the Sobolev-Troisi inequality \cite{tro}
\begin{lem}
    Let $\Omega\subset \mathbb{R}^{N}$ be a bounded open set and consider $u\in W^{1,p_{i}}(\Omega)$, $p_{i}\geq1$ for all $i=1,..,N.$ Set
    \begin{equation*}
       \frac{1}{\Bar{p}}=\frac{1}{N}\sum_{i=1}^{N}\frac{1}{p_{i}},~\Bar{p}^{*}=\frac{N\Bar{p}}{N-\Bar{p}}. 
    \end{equation*}
If $\Bar{p}<N$, there exists a positive constant $C$ depending only on $\Omega,~p_{i}$ and $N$ such that
\begin{equation}
    \|u\|_{\Bar{p}^{*}}\leq C\prod_{i=1}^{N}\left\|\frac{\partial u}{\partial x_{i}} \right\|^{\frac{1}{N}}_{L^{p_{i}}(\Omega)}.
\end{equation}
\end{lem}
Next, we state the following anisotropic embedding which can be found in \cite{yu}.
\begin{lem}
    Let $u\in C(0,T; L^{2}(\Omega))\cap L^{\Vec{p}}(0,T;W^{1,\Vec{p}}_{0}(\Omega))$ and assume that
    \begin{equation*}
        \frac{2N}{N+2}\leq\Bar{p}<N,~l=\Bar{p}(1+\frac{2}{N}).
    \end{equation*}
    Then, there exists a constant $C>0$ such that
    \begin{equation}
        \int_{\Omega_{T}}|u|^{l}~dxdt\leq C\left(\underset{t\in[0,T]}{\sup}\int_{\Omega}|u|^{2}~dx+\sum_{i=1}^{N}\int_{\Omega_{T}}\left|\frac{\partial u}{\partial x_{i}} \right|^{p_{i}}~dxdt
        \right)^{\frac{N+\Bar{p}}{N}}.
    \end{equation}
\end{lem}

To formally define the weak solutions to obstacle problems related to (1.1), we consider the following class of functions
\begin{equation*}
    K_{\phi}(\Omega):=\{u\in C(0,T; L^{2}(\Omega));~u\in L^{\Vec{p}}(0,T;~W^{1,\Vec{p}}(\Omega)),~u\geq\phi~\text{a.e. in}~\Omega_{T}\}.
\end{equation*}
Furthermore, the class of admissible comparison functions is defined  as follows 
\begin{equation*}
    K'_{\phi}(\Omega):=\{u\in K_{\phi}(\Omega),~u_{t}\in L^{2}(\Omega_{T})    \}.
\end{equation*}
Now, we have enough tools to give a definition of weak solutions to our obstacle problem. The existence of such solutions is guaranteed by \cite{bog3,ham3}.
\begin{defn} 
    We define $u\in K_{\phi}(\Omega_{T})$ as a local weak solution to the obstacle problems associated with (1.1) if and only if
\begin{equation}
\begin{split}
    <<u_{t},\varphi(v-u)>>+\sum_{i=1}^{N}\int_{\Omega}\left|\frac{\partial u}{\partial x_{i}} \right|^{p_{i}-2} \frac{\partial u}{\partial x_{i}}\frac{\partial}{\partial x_{i}}(\varphi(v-u))~dxdt
    \geq 0
     \end{split}
     \end{equation}
holds true for all comparison functions $v\in K'_{\phi}(\Omega_{T})$ and every test function $\varphi\in C_{0}^{\infty}(\Omega_T, \mathbb{R}^{+})$. The time term above is defined as
\begin{equation*}
   <<u_{t},\varphi(v-u)>>= \int_{\Omega_{T}}\left\{\varphi_{t}\left[\frac{1}{2}u^{2}-uv \right]-\varphi uv_{t}   \right\}~dxdt  
\end{equation*}
\end{defn}
To address the potential lack of differentiability in time of weak solutions to obstacle problems related to (1.1), the following time mollification has been proven to be useful
\begin{equation}
 [u]_{h}(x,t):=\frac{1}{h}\int_{0}^{t}e^{\frac{s-t}{h}}u(x,s)~ds,  
\end{equation}
for $u\in L^{1}(\Omega_{T})$ and $h>0$. We summarize some elementary properties of (2.4), which can be retrieved from \cite{kin} in the following lemma.

\begin{lem}
    For $p\geq1$, we have
    \begin{itemize}
        \item If $u\in L^{p}(\Omega_{T})$ then $[u]_{h}\longrightarrow u$ in $L^{p}(\Omega_{T})$ as $h\downarrow0$ and
        \begin{equation*}
            \partial_{t}[u]_{h}=\frac{1}{h}(u-[u]_{h})\in L^{p}(\Omega) ~\text{for every } h>0.  
        \end{equation*} 
        \item If $\nabla u \in L^{p}(\Omega_{T}, \mathbb{R}^{N})$ then $\nabla[u]_{h}=[\nabla u]_{h}$ and $\nabla[u]_{h}\longrightarrow \nabla u$ in $ L^{p}(\Omega_{T}, \mathbb{R}^{N})$ as $h\downarrow0$.
        \item If $u\in C^{0}(\overline{\Omega_{T}})$, then $[u]_{h}\longrightarrow u$ uniformly in $\Omega_{T}$ as $h\downarrow0$.
        
    \end{itemize}
\end{lem}
Finally, we present the following technical lemma which is frequently used in this work
\begin{lem}
    Let $(X_{i})_{i\in \mathbb{N}}$ be a sequence of positive real numbers with
    \begin{equation*}
        X_{i+1}\leq C B^{i} X_{i}^{1+\alpha},~\text{for all}~i\in \mathbb{N},
    \end{equation*}
    for constants $C,\alpha>0$ and $B>1$. Then
    \begin{equation*}
        X_{0}\leq C^{-\frac{1}{\alpha}}B^{-\frac{1}{\alpha^{2}}}
    \end{equation*}
    implies $X_{i}\longrightarrow0$ as $i\longrightarrow\infty$.
\end{lem}
\section{Local energy and logarithmic estimates}

Customarily, we use the symbols $(u-k)_+$ and $(u-k)_-$ to denote the positive and negative parts of these truncated functions, respectively, such that for $k>0$
\begin{equation*}
    (u-k)_+=\max\{u-k,0\}~~~\text{and}~~~(u-k)_{-}=\max\{k-u,0\}.
\end{equation*}
To make the calculations easier, we will consider cylinders with the vertex at the origin $(0,0)$. The results for cylinders with a vertex at a different point $(x_0,t_0)$ can be obtained by simply translating the calculations.

For $\rho>0$, define the anisotropic box and cylinder:
\[
K_\rho := \prod_{i=1}^N (-L_i(\rho), L_i(\rho)), \qquad 
Q(s,\rho) := K_\rho \times (-s, 0),
\]
where $L_i(\rho)$ denotes the half-length of $K_\rho$ in the $x_i$-direction.

\begin{lem}
Let $Q(s,\rho) \subset \Omega_T$ be a cylinder as above, and let $\phi\in C^{0}(\Omega_T)$. Let $u\in K_{\phi}(\Omega_T)$ be a weak local solution of (1.1) in the sense of Definition 2.3. For any function $\xi \in C^\infty(Q(s,\rho); \mathbb{R}_{\ge 0})$ vanishing on the lateral boundary $\partial K_\rho \times (-s,0)$, there exists a constant $C>0$ depending on the data such that for any $\alpha \geq p^+$, the following estimates hold:
\begin{enumerate}
    \item For any $k>0$, we have
    \begin{equation}
\begin{split}
\underset{-s<t<0}{\sup}\int_{K_{\rho}}&\xi^{\alpha}(u-k)_{-}^{2}~dx+C\sum_{i=1}^{N}\int_{Q(s,\rho)}\xi^{\alpha}\left|\frac{\partial}{\partial x_{i}}(u-k)_{-} \right|^{p_{i}}~dxdt\\
&\leq\int_{K_{\rho}\times\{-s\}}\xi^{\alpha}(u-k)_{-}^{2}~dx+C\sum_{i=1}^{N}\biggl\{\int_{Q(s,\rho)} \left|\frac{\partial\xi}{\partial t} \right|(u-k)_{-}^{2}~dxdt\\
&+\int_{Q(s,\rho)} \left|\frac{\partial\xi}{\partial x_{i}} \right|^{p_{i}}(u-k)_{-}^{p_{i}}~dxdt\biggr\}.
\end{split}
\label{eq:energy-neg}
\end{equation}
    \item For all $k\geq\underset{Q(s,\rho)}{\sup}~\phi$, we have
\begin{equation}
\begin{split}
\underset{-s<t<0}{\sup}\int_{K_{\rho}}&\xi^{\alpha}(u-k)_{+}^{2}~dx+C\sum_{i=1}^{N}\int_{Q(s,\rho)}\xi^{\alpha}\left|\frac{\partial}{\partial x_{i}}(u-k)_{+} \right|^{p_{i}}~dxdt\\
&\leq\int_{K_{\rho}\times\{-s\}}\xi^{\alpha}(u-k)_{+}^{2}~dx+C\sum_{i=1}^{N}\biggl\{\int_{Q(s,\rho)} \left|\frac{\partial\xi}{\partial t} \right|(u-k)_{+}^{2}~dxdt\\
&+\int_{Q(s,\rho)} \left|\frac{\partial\xi}{\partial x_{i}} \right|^{p_{i}}(u-k)_{+}^{p_{i}}~dxdt
\biggr\}.
\end{split}
\label{eq:energy-pos}
\end{equation}
\end{enumerate}
\end{lem}

\begin{proof}
We begin by introducing two nonnegative piecewise smooth functions. For the space-time cutoff function, we take a function of the product form
\[
\xi(x,t) = \eta(t) \prod_{i=1}^N \xi_i(x_i),
\]
where each $\xi_i \in C_0^\infty((-L_i(\rho), L_i(\rho)), \mathbb{R}^+)$ vanishes on $\partial(-L_i(\rho), L_i(\rho))$ (so that $\xi$ vanishes on the lateral boundary $\partial K_\rho \times (-s,0)$), and $\eta \in C^1([-s,0], \mathbb{R}^+)$. The time localization is handled separately by $\psi_{\varepsilon} \in W_{0}^{1,\infty}((-s,0);~[0,1])$, which satisfies
\begin{equation*}
  \psi_{\varepsilon}(t)= \begin{cases}
     0,~~~~~~~~~~~&~~\text{for}~-s\leq t\leq s_{1}-\varepsilon,\\
     1+\frac{t-s_1}{\varepsilon},~~~~~~~~~~~&~~\text{for}~s_{1}-\varepsilon< t\leq s_{1},\\
     1,~~~~~~~~~~~&~~\text{for}~s_{1}< t< s_{2},\\
     1-\frac{t-s_2}{\varepsilon},~~~~~~~~~~~&~~\text{for}~s_{2}\leq t< s_{2}+\varepsilon,\\
     0,~~~~~~~~~~~&~~\text{for}~s_{2}+\varepsilon\leq t\leq0.\\
  \end{cases}  
\end{equation*}
Furthermore, by using the time mollification $[.]_{h}$ defined in (2.4), for $k>0$, and $h>0$, let
\begin{equation}
    v_{h}=[u]_{h}+([u]_{h}-k)_{-}+\left\|\phi-[\phi]_{h} \right\|_{L^{\infty}(Q(s,\rho))}.
\end{equation}
From (3.3), we deduce that $v_{h}\ge \phi$ a.e. in $Q(s,\rho)$. Indeed, since $([u]_h - k)_- \ge 0$, the positivity-preserving property $[u]_h \ge [\phi]_h$ a.e. together with $[\phi]_h + \|\phi - [\phi]_h\|_{L^\infty(Q(s,\rho))} \ge \phi$ a.e. yields
\begin{equation*}
    v_{h} \ge [u]_h + \|\phi - [\phi]_h\|_{L^\infty(Q(s,\rho))} \ge [\phi]_h + \|\phi - [\phi]_h\|_{L^\infty(Q(s,\rho))} \ge \phi \quad \text{a.e. in } Q(s,\rho).
\end{equation*}
Therefore, taking $v_{h}$ as an admissible comparison function and choosing the cut-off function $\varphi = \xi^{\alpha}(\psi_{\varepsilon})_{\delta}$ in the variational inequality (2.3), we obtain
\begin{equation}
\begin{aligned}
    \int_{Q(s,\rho)} \partial_{t}\left(\xi^{\alpha}(\psi_{\varepsilon})_{\delta}\right) &\left(\frac{1}{2}u^{2}-uv_{h}\right) dx dt - \int_{Q(s,\rho)} \xi^{\alpha}(\psi_{\varepsilon})_{\delta} \, u \, \partial_{t}v_{h} \, dx dt \\
    &+ \sum_{i=1}^{N}\int_{Q(s,\rho)}\left|\frac{\partial u}{\partial x_{i}} \right|^{p_{i}-2}\frac{\partial u}{\partial x_{i}}(\psi_{\varepsilon})_{\delta} \frac{\partial}{\partial x_{i}}\left( \xi^{\alpha} (v_{h}-u) \right) dx dt \ge 0,
\end{aligned}
\end{equation}
where $\alpha > 0$ is a constant to be determined later, and $(\psi_{\varepsilon})_{\delta}$ is a mollification of $\psi_{\varepsilon}$ defined in [\cite{giu}, Section 2.2] with $0 < \delta < \frac{\varepsilon}{2}$, and $\mathrm{supp}\left(\xi^{\alpha}(\psi_{\varepsilon})_{\delta}\right) \subset Q(s,\rho)$. Next, in order to simplify the second integral on the left-hand side of (3.4), we use the following assertion which we obtained from Lemma 2.4
\begin{equation}
    \partial_{t}v_{h}= \frac{\partial}{\partial t}\left( [u]_{h} + ([u]_{h} - k)_{-} \right)=\begin{cases}
        \frac{1}{h}(u-[u]_{h})&~~~\text{if}~Q(s,\rho)\cap\{[u]_{h}>k\},\\
        0&~~~\text{otherwise}.
       \end{cases}
\end{equation}
As a result, we obtain
\begin{equation}
\begin{aligned}
    \int_{Q(s,\rho)} \xi^{\alpha}(\psi_{\varepsilon})_{\delta} \, u \, \partial_{t}v_{h} \, dx dt 
    &= \int_{Q(s,\rho)} \xi^{\alpha}(\psi_{\varepsilon})_{\delta} (u - [u]_{h}) \, \partial_{t}v_{h} \, dx dt \\
    &\quad + \int_{Q(s,\rho)} \xi^{\alpha}(\psi_{\varepsilon})_{\delta} [u]_{h} \, \partial_{t}v_{h} \, dx dt \\
    &\ge \int_{Q(s,\rho)} \xi^{\alpha}(\psi_{\varepsilon})_{\delta} [u]_{h} \, \frac{\partial}{\partial t}\left( [u]_{h} + ([u]_{h} - k)_{-} \right) dx dt,
\end{aligned}
\end{equation}
where we used the fact that
\begin{equation*}
    \int_{Q(s,\rho)} \xi^{\alpha}(\psi_{\varepsilon})_{\delta} (u - [u]_{h}) \, \partial_{t}v_{h} \, dx dt = \frac{1}{h} \int_{Q(s,\rho) \cap \{[u]_{h} > k\}} \xi^{\alpha}(\psi_{\varepsilon})_{\delta} (u - [u]_{h})^{2} \, dx dt \ge 0.
\end{equation*}
Afterward, the last term on the right-hand side of (3.6) becomes
\begin{equation*}
\begin{split}
\int_{Q(s,\rho)}&\xi^{\alpha}(\psi_{\varepsilon})_{\delta}[u]_{h}\frac{\partial}{\partial t}([u]_{h}+([u]_{h}-k)_{-})~dxdt\\
=&-\int_{Q(s,\rho)}\frac{\partial}{\partial t}(\xi^{\alpha}(\psi_{\varepsilon})_{\delta})[u]_{h}^{2}~dxdt-
\int_{Q(s,\rho)}\frac{\partial}{\partial t}(\xi^{\alpha}(\psi_{\varepsilon})_{\delta})[u]_{h}([u]_{h}-k)_{-}~dxdt\\
&-\frac{1}{2}\int_{Q(s,\rho)}\xi^{\alpha}(\psi_{\varepsilon})_{\delta}\frac{\partial[u]^{2}_{h}}{\partial t}~dxdt
-\int_{Q(s,\rho)}\xi^{\alpha}(\psi_{\varepsilon})_{\delta}\frac{\partial [u]_{h}}{\partial t}([u]_{h}-k)_{-}~dxdt\\
=&-\frac{1}{2}\int_{Q(s,\rho)}\frac{\partial}{\partial t}(\xi^{\alpha}(\psi_{\varepsilon})_{\delta})[u]_{h}^{2}~dxdt
-\int_{Q(s,\rho)}\frac{\partial}{\partial t}(\xi^{\alpha}(\psi_{\varepsilon})_{\delta})[u]_{h}([u]_{h}-k)_{-}~dxdt\\
&+\int_{Q(s,\rho)}\xi^{\alpha}(\psi_{\varepsilon})_{\delta}\frac{\partial}{\partial t}\int_{0}^{([u]_{h}-k)_{-}}\tau~d\tau dxdt\\
=&-\int_{Q(s,\rho)}\frac{\partial}{\partial t}(\xi^{\alpha}(\psi_{\varepsilon})_{\delta})(\frac{1}{2}[u]_{h}^{2}+[u]_{h}([u]_{h}-k)_{-})~dxdt
-\int_{Q(s,\rho)}\frac{\partial}{\partial t}(\xi^{\alpha}(\psi_{\varepsilon})_{\delta})\int_{0}^{([u]_{h}-k)_{-}}\tau~d\tau dxdt.
\end{split}
\end{equation*}
In conclusion, (3.6) becomes
\begin{equation}
    \begin{split}
      \int_{Q(s,\rho)} \xi^{\alpha}(\psi_{\varepsilon})_{\delta} u\partial_{t} v_{h}~dxdt
     \geq&-\int_{Q(s,\rho)}\frac{\partial}{\partial t}(\xi^{\alpha}(\psi_{\varepsilon})_{\delta})(\frac{1}{2}[u]_{h}^{2}+[u]_{h}([u]_{h}-k)_{-})~dxdt\\
&-\int_{Q(s,\rho)}\frac{\partial}{\partial t}(\xi^{\alpha}(\psi_{\varepsilon})_{\delta})\int_{0}^{([u]_{h}-k)_{-}}\tau~d\tau dxdt.
\end{split}
\end{equation}
Therefore, passing to the limit as $h \downarrow 0$ and noting that $v_h \to u + (u-k)_{-}$ pointwise almost everywhere in $Q(s,\rho)$, the sum of the first two time-derivative terms in the left hand side of (3.4) simplifies to:
\begin{equation}
\begin{aligned}
    \lim_{h\downarrow 0} &\int_{Q(s,\rho)} \frac{\partial}{\partial t}\left( \xi^{\alpha}(\psi_{\varepsilon})_{\delta} \right) \left( \frac{1}{2}u^{2} - u v_{h} \right) dx dt - \int_{Q(s,\rho)} \xi^{\alpha}(\psi_{\varepsilon})_{\delta} \, u \, \partial_{t} v_{h} \, dx dt \\
    &\le \int_{Q(s,\rho)} \frac{\partial}{\partial t}\left( \xi^{\alpha}(\psi_{\varepsilon})_{\delta} \right) \left( -\frac{1}{2}u^2 - u(u-k)_{-} + \frac{1}{2}u^2 + u(u-k)_{-} + \int_{0}^{(u-k)_{-}} \tau \, d\tau \right) dx dt \\
    &= \int_{Q(s,\rho)} \frac{\partial}{\partial t}\left( \xi^{\alpha}(\psi_{\varepsilon})_{\delta} \right) \int_{0}^{(u-k)_{-}} \tau \, d\tau \, dx dt\\
    &=\int_{Q(s,\rho)} \frac{\partial}{\partial t}\left( \xi^\alpha (\psi_\varepsilon)_\delta \right) \frac{1}{2} \left((u-k)_-\right)^2\, dx dt. 
\end{aligned}
\end{equation}
Next, we simplify the spatial derivative integral as $h \downarrow 0$. By Lemma 2.4, we have the strong convergence
\begin{equation*}
    \frac{\partial}{\partial x_{i}}\big(\xi^{\alpha}(v_{h}-u)\big) \underset{h\downarrow 0}{\longrightarrow} \frac{\partial}{\partial x_{i}}\big(\xi^{\alpha}(u-k)_{-}\big) \quad \text{in } L^{p_{i}}(Q(s,\rho))
\end{equation*}
for each $i = 1, \dots, N$. Passing to the limit, expanding the derivative $\partial_{x_i}\big(\xi^\alpha (u-k)_-\big)$, and using the relation $\partial_{x_i} u = -\partial_{x_i} (u-k)_{-}$ on $\{u < k\}$, we obtain
\begin{equation*}
\begin{aligned}
    \lim_{h \downarrow 0} &\int_{Q(s,\rho)} (\psi_{\varepsilon})_{\delta} \left|\frac{\partial u}{\partial x_{i}}\right|^{p_{i}-2} \frac{\partial u}{\partial x_{i}} \frac{\partial}{\partial x_{i}} \big(\xi^{\alpha}(v_{h}-u)\big) \, dx dt \\
    &= \int_{Q(s,\rho)} (\psi_{\varepsilon})_{\delta} \left|\frac{\partial u}{\partial x_{i}}\right|^{p_{i}-2} \frac{\partial u}{\partial x_{i}} \frac{\partial}{\partial x_{i}} \big(\xi^{\alpha}(u-k)_{-}\big) \, dx dt \\
    &= -\int_{Q(s,\rho)} (\psi_{\varepsilon})_{\delta} \xi^{\alpha} \left|\frac{\partial}{\partial x_{i}}(u-k)_{-}\right|^{p_{i}} dx dt \\
    &\quad - \alpha \int_{Q(s,\rho)} (\psi_{\varepsilon})_{\delta} \xi^{\alpha-1} (u-k)_{-} \left|\frac{\partial}{\partial x_{i}}(u-k)_{-}\right|^{p_{i}-2} \frac{\partial (u-k)_-}{\partial x_i} \frac{\partial \xi}{\partial x_i} \, dx dt.
\end{aligned}
\end{equation*}
Applying Young's inequality with parameter $\nu > 0$ to the cross term yields
\begin{equation*}
\begin{aligned}
    \alpha \int_{Q(s,\rho)} &(\psi_{\varepsilon})_{\delta} \xi^{\alpha-1} (u-k)_{-} \left|\frac{\partial}{\partial x_{i}}(u-k)_{-}\right|^{p_{i}-1} \left|\frac{\partial \xi}{\partial x_i}\right| dx dt \\
    &\le \alpha\nu \int_{Q(s,\rho)} (\psi_{\varepsilon})_{\delta} \xi^{(\alpha-1)\frac{p_{i}}{p_{i}-1}} \left|\frac{\partial}{\partial x_{i}}(u-k)_{-}\right|^{p_{i}} dx dt \\
    &\quad + \alpha C(\nu) \int_{Q(s,\rho)} (\psi_{\varepsilon})_{\delta} \left|\frac{\partial \xi}{\partial x_{i}}\right|^{p_{i}} (u-k)_{-}^{p_{i}} \, dx dt.
\end{aligned}
\end{equation*}
Since $\alpha > p_+ \ge p_i$, we have $(\alpha-1)\frac{p_i}{p_i-1} \ge \alpha$. Because $0 < \xi \le 1$, it follows that $\xi^{(\alpha-1)\frac{p_{i}}{p_{i}-1}} \le \xi^{\alpha}$. Thus, choosing $\nu > 0$ small enough such that $\alpha\nu < 1$, we arrive at the bound
\begin{equation}
\begin{aligned}
    \lim_{h \downarrow 0} &\int_{Q(s,\rho)} (\psi_{\varepsilon})_{\delta} \left|\frac{\partial u}{\partial x_{i}}\right|^{p_{i}-2} \frac{\partial u}{\partial x_{i}} \frac{\partial}{\partial x_{i}} \big(\xi^{\alpha}(v_{h}-u)\big) \, dx dt \\
    &\le -C \int_{Q(s,\rho)} (\psi_{\varepsilon})_{\delta} \xi^{\alpha} \left|\frac{\partial}{\partial x_{i}}(u-k)_{-}\right|^{p_{i}} dx dt + \alpha C(\nu) \int_{Q(s,\rho)} (\psi_{\varepsilon})_{\delta} \left|\frac{\partial \xi}{\partial x_{i}}\right|^{p_{i}} (u-k)_{-}^{p_{i}} \, dx dt,
\end{aligned}
\end{equation}
where $C = 1 - \alpha\nu > 0$. Therefore, by putting (3.8), and (3.9) into (3.4) we get
\begin{equation}
    \begin{split}
\sum_{i=1}^{N}&\int_{Q(s,\rho)}(\psi_{\varepsilon})_{\delta}\xi^{\alpha}\left|\frac{\partial}{\partial x_{i}}(u-k)_{-} \right|^{p_{i}}~dxdt\leq
\int_{Q(s,\rho)}\frac{\partial}{\partial t}((\psi_{\varepsilon})_{\delta}\xi^{\alpha})\int_{0}^{(u-k)_{-}}\tau d\tau dxdt\\
&+C\sum_{i=1}^{N}\int_{Q(s,\rho)}(\psi_{\varepsilon})_{\delta}\left|\frac{\partial\xi}{\partial x_{i}} \right|^{p_{i}}(u-k)_{-}^{p_{i}}~dxdt
\end{split}
\end{equation}
Fix $-s < t_1 < t_2 < 0$. Setting $s_1 = t_1$ and $s_2 = t_2$ in the definition of $\psi_\varepsilon(t)$, we evaluate the double limit $\delta \to 0$ and $\varepsilon \to 0$ on the time-derivative term:
\begin{equation}
\begin{aligned}
\lim_{\varepsilon \to 0} \lim_{\delta \to 0} &\int_{Q(s,\rho)} \frac{\partial}{\partial t}\left( \xi^\alpha (\psi_\varepsilon)_\delta \right) \frac{1}{2} \left((u-k)_-\right)^2 dx dt \\
&= \lim_{\varepsilon \to 0} \int_{-s}^0 \partial_t \psi_\varepsilon(t) \, G(t) \, dt + \frac{\alpha}{2} \int_{t_1}^{t_2} \int_{K_\rho} \xi^{\alpha-1} \frac{\partial \xi}{\partial t} \left((u-k)_-\right)^2 dx dt \\
&= G(t_1) - G(t_2) + \frac{\alpha}{2} \int_{t_1}^{t_2} \int_{K_\rho} \xi^{\alpha-1} \frac{\partial \xi}{\partial t} \left((u-k)_-\right)^2 dx dt,
\end{aligned}
\end{equation}
where $G(t) := \frac{1}{2} \int_{K_\rho} \xi^\alpha(x,t) \left((u(x,t)-k)_-\right)^2 dx$, and as $\delta,\varepsilon\to0$, the mollified cutoff $(\psi_\varepsilon)_\delta$ converges almost everywhere to $\chi_{(t_1,t_2)}$, while
\[
\partial_t\psi_\varepsilon\xrightarrow{\mathcal{D}'}\delta_{t_1}-\delta_{t_2}.
\]
This extracts the temporal boundary values:
\begin{equation*}
\lim_{\varepsilon \to 0} \int_{-s}^{0} \partial_t \psi_\varepsilon(t) G(t) \, dt = \langle \delta_{t_1} - \delta_{t_2}, \, G \rangle = G(t_1) - G(t_2).
\end{equation*}
Applying these limits to the full weak formulation yields the intermediate energy estimate over $(t_1, t_2)$:
\begin{equation*}
\begin{aligned}
G(t_2) + C \sum_{i=1}^{N} \int_{t_1}^{t_2} \int_{K_\rho}& \xi^{\alpha} \left| \frac{\partial (u-k)_{-}}{\partial x_{i}} \right|^{p_{i}} dx dt \\
\le& G(t_1) + C \sum_{i=1}^{N} \int_{t_1}^{t_2} \int_{K_\rho} \left| \frac{\partial \xi}{\partial x_{i}} \right|^{p_{i}} (u-k)_{-}^{p_{i}} dx dt\\
&+ \frac{\alpha}{2} \int_{t_1}^{t_2} \int_{K_\rho} \xi^{\alpha-1} \left| \frac{\partial \xi}{\partial t} \right| \left((u-k)_-\right)^2 dx dt.
\end{aligned}
\end{equation*}
Passing to the limit $t_1 \downarrow -s$, the spatial energy term $G(t_1)$ converges by continuity of the trace to the initial slice value $G(-s) = \frac{1}{2}\int_{K_{\rho}\times\{-s\}}\xi^{\alpha}(u-k)_{-}^{2}~dx$. This extends the lower integration limit in the space-time gradient integral down to $-s$. Next, taking the supremum over $t_2 \in (-s, 0)$ bounds the spatial norm term, while letting $t_2 \uparrow 0$  covering the full domain $Q(s,\rho)$, and thus, we get the desired (3.1).

The proof of (3.2) is similar to the one of (3.1) where we took $$v_{h}=[u]_{h}-([u]_{h}-k)_{+}+\|\phi-[\phi]_{h}\|_{L^{\infty} (Q(s,\rho))}$$ as an admissible test function in (2.3) for any fixed $k\geq \underset{Q(s,\rho)}{\sup}\phi$, replacing $([u]_{h}-k)_{-}$ by $([u]_{h}-k)_{+}$ and using the following identity
\begin{equation*}
    \partial_{t}[u]_{h}([u]_{h}-k)_{+}=\frac{\partial}{\partial t}\int_{0}^{([u]_{h}-k)_{+}}\tau d\tau.
\end{equation*}
Hence we get the desired result.
\end{proof}

Now, we will introduce the following logarithmic function
\begin{equation}
    \Gamma_{\pm}(u)=\Gamma(H^{\pm}_{k},(u-k)_{\pm},c)=\left[\ln\left( \frac{H^{\pm}_{k}}{(H^{\pm}_{k}-(u-k)_{\pm}+c} \right) \right]_{+}
\end{equation}
where $H^{\pm}_{k}=\underset{Q(s,\rho)}{ess\sup}|(u-k)_{\pm}|$ and $0<c<H^{\pm}_{k}$. Then, we get the following properties
\begin{equation}
    \begin{cases}
        \Gamma_{\pm}(u)=0,&~~~~\text{if}~(u-k)_{\pm}\leq c,\\
        0\leq \Gamma_{\pm}(u)\leq \ln\left(\frac{H^{\pm}_{k}}{c} \right),&~~~~\text{if}~(u-k)_{\pm}\leq H^{\pm}_{k},\\
        0\leq \Gamma'_{\pm}(u)\leq \frac{1}{c},&~~~~\text{if}~(u-k)_{\pm}\neq c,\\
        \left(\Gamma_{\pm}(u)\right)^{''}=\left(\Gamma'_{\pm}(u) \right)^{2},&~~~~\text{if}~(u-k)_{\pm}\neq c.
    \end{cases}
\end{equation}
Moreover, since $\Gamma_{\pm}^{2}$ is differentiable in $[0,T]$ we have the following
\begin{equation}
    \left(\Gamma_{\pm}^{2} \right)'=2\Gamma_{\pm}\Gamma'_{\pm}~\text{and}~\left(\Gamma_{\pm}^{2}\right)^{''}=2(1+\Gamma_{\pm})\left(\Gamma'_{\pm}\right)^{2} ~\text{in}~[0,H_{k}^{\pm}].
\end{equation}
Also, by using Theorem 4.1, we may assume that
\begin{equation}
   H^{\pm}_{k}=\underset{Q(s,\rho)}{ess\sup}|(u-k)_{\pm}|<\infty. 
\end{equation}

\begin{lem}
Let $K_{\rho_2} \subset K_{\rho_1}$ be nested anisotropic boxes with half-lengths $L_i(\rho_1)$ and $L_i(\rho_2)$, respectively. Let $0<t_1<t_2<T$. We abbreviate $Q_1 := K_{\rho_1} \times (t_1,t_2)$. For $H^{\pm}_{k}=\underset{Q_1}{ess\sup}|(u-k)_{\pm}|$, and $\phi\in C^{0}(\Omega_T)$, let $u\in K_{\phi}(\Omega_T)$ be a weak local solution of (1.1) in the sense of Definition 2.3. Then, there exists a constant $C>0$ depending on the data such that the following estimates hold
\begin{enumerate}
    \item for $k\geq\underset{Q_1}{\sup}~\phi$ we have
        \begin{equation}
    \begin{split}
\underset{t\in(t_1,t_2)}{ess\sup}\int_{K_{\rho_2}}\Gamma_{+}^{2}~dx\leq \int_{K_{\rho_1}\times\{t_1\}}\Gamma^{2}_{+}~dx+
 C\sum_{i=1}^{N}\frac{1}{(L_i(\rho_1)-L_i(\rho_2))^{p_i}}\int_{Q_1}\Gamma_{+}\left(\Gamma'_{+}\right)^{2-p_i}~dxdt.
 \end{split}
\label{eq:log-pos}
    \end{equation}
     \item For any $k>0$, we have
       \begin{equation}
        \begin{split}
\underset{t\in(t_1,t_2)}{ess\sup}&\int_{K_{\rho_2}}\Gamma_{-}^{2}~dx\leq \int_{K_{\rho_1}\times\{t_1\}}\Gamma^{2}_{-}~dx
+ C\sum_{i=1}^{N}\frac{1}{(L_i(\rho_1)-L_i(\rho_2))^{p_i}}\int_{Q_1}\Gamma_{-}\left(\Gamma'_{-}\right)^{2-p_i}~dxdt.
 \end{split}
\label{eq:log-neg}
    \end{equation}
    \end{enumerate}
\end{lem}

\begin{proof}
We begin the proof by letting
\begin{equation}
    v_{h}=[u]_{h}-\beta\left(\Gamma_{+}^{2} \right)'([u]_{h})+\left\|\phi-[\phi]_{h}  \right\|_{L^{\infty}(\Omega_{T})},
\end{equation}
where $\beta>0$ is a constant to be specified later. For $k\geq\underset{Q_{1}}{\sup}\phi$, we have that $v_{h}\geq\phi$. Indeed, if $[u]_{h}<k$, we get that
\begin{equation}
    v_{h}=[u]_{h}+\left\|\phi-[\phi]_{h}\right\|_{L^{\infty}(Q_{1})}\geq \left\|\phi\right\|_{L^{\infty}(Q_{1})}\geq \phi~\text{a.e. in}~Q_{1}.
\end{equation}
Else, for $0<\beta\leq\frac{k-\underset{Q_{1}}{\sup \phi}}{\underset{[0,H^{+}_{k}]}{\sup}\left(\Gamma_{+}^{2} \right)'}$, we get
\begin{equation}
    v_{h}\geq k-\beta\left(\Gamma_{+}^{2} \right)'([u]_{h})\geq k-\beta\underset{[0,H^{+}_{k}]}{\sup}\left(\Gamma_{+}^{2} \right)'\geq \underset{Q_{1}}{\sup}\phi\geq\phi~\text{a.e. in}~Q_{1}.
\end{equation}
Moreover, since $\left(\Gamma_{+}^{2} \right)'$ is a Lipschitz function and be Lemma 2.4, we can use $v_{h}$ as an admissible comparison function in (2.3) such that
\begin{equation}
    \partial_{t}v_{h}=\frac{1}{h}(u-[u]_{h})\left(1-\beta\left(\Gamma_{+}^{2} \right)^{''}([u]_{h})\right),
\end{equation}
where the terms involving $\left(\Gamma_{+}^{2} \right)^{''}$ is well defined since $\partial_{t}[u]_{h}=0$ in the set $\{([u]_{h}-k)_{+}=~C\in\mathbb{R}\}$. 
Thereafter, we are going to simplify the second integral on the left-hand side of (2.3). Let $\xi\in C^{1}_{0}(K_{\rho_1},\mathbb{R}^+)$ be a spatial cutoff of the product form $\xi(x)=\prod_{i=1}^N \xi_i(x_i)$ with $\xi=1$ on $K_{\rho_2}$ and
\[
\left|\frac{\partial\xi}{\partial x_{i}}\right|\leq \frac{C}{L_i(\rho_1)-L_i(\rho_2)} \quad \text{for all } i=1,..,N.
\]
For the time direction, we take a cutoff $\psi_{\varepsilon}$ adapted as in the proof of  Lemma 3.1. Then we have
\begin{equation}\begin{split}
    \int_{Q_1}&\xi^{\alpha}\psi_{\varepsilon}u\partial_{t}v_{h}~dxdt\\
    &= \int_{Q_1}\xi^{\alpha}\psi_{\varepsilon}(u-[u]_{h})\partial_{t}v_{h}~dxdt
    + \int_{Q_1}\xi^{\alpha}\psi_{\varepsilon}[u]_{h}\partial_{t} v_{h}~dxdt.
 \end{split}
\end{equation}
Then, for $\beta\leq\frac{1}{\underset{[0,H^{+}_{k}]}{\sup}\left(\Gamma_{+}^{2}\right)^{''} }$, the first integral on the right-hand side of (3.22) becomes
\begin{equation}
  \frac{1}{h}\int_{Q_1}\xi^{\alpha}\psi_{\varepsilon}(u-[u]_{h})(u-[u]_{h})\left(1-\beta\left( \Gamma_{+}^{2}\right)^{''}([u]_{h}) \right)~dxdt\geq0.
\end{equation}
Next, for the second integral on the right-hand side of (3.22), we have
\begin{equation}
    \begin{split}
        \int_{Q_1}&\xi^{\alpha}\psi_{\varepsilon}[u]_{h}\partial_{t}v_{h}~dxdt = \frac{1}{2}\int_{Q_1}\xi^{\alpha}\psi_{\varepsilon}\partial_{t}[u]_{h}^{2}~dxdt + \int_{Q_1}\xi^{\alpha}\psi'_{\varepsilon}[u]_{h}\beta\left(\Gamma_{+}^{2} \right)'([u]_{h})~dxdt \\
        &+\int_{Q_1}\xi^{\alpha}\psi_{\varepsilon}\partial_{t}[u]_{h}\beta\left(\Gamma_{+}^{2} \right)'([u]_{h})~dxdt \\
        =& -\frac{1}{2}\int_{Q_1}\xi^{\alpha}\psi'_{\varepsilon}[u]_{h}^{2}~dxdt + \int_{Q_1}\xi^{\alpha}\psi'_{\varepsilon}[u]_{h}\beta\left(\Gamma_{+}^{2} \right)'([u]_{h})~dxdt \\
        &+\int_{Q_1}\xi^{\alpha}\psi_{\varepsilon}\partial_{t}\left( \int_{0}^{([u]_{h}-k)_{+}}\beta\left( \Gamma_{+}^{2}\right)'(s)~ds\right)~dxdt \\
        =& \int_{Q_1}\xi^{\alpha}\psi'_{\varepsilon}\left( [u]_{h}\beta\left(\Gamma_{+}^{2} \right)'([u]_{h})-\frac{1}{2}[u]_{h}^{2}\right)~dxdt \\
        &-\int_{Q_1}\xi^{\alpha}\psi'_{\varepsilon}\int_{0}^{([u]_{h}-k)_{+}}\beta\left( \Gamma_{+}^{2}\right)'(s)~dsdxdt.
    \end{split}
\end{equation}
Therefore, by putting (3.24) into (3.22), the first term on the left-hand side of (2.3) becomes (after integration by parts and letting $h\to0$, then $\varepsilon\to0$)
\begin{equation}
    \begin{split}
        \underset{h\downarrow0}{\lim}\int_{Q_1}&\xi^{\alpha}\psi'_{\varepsilon} \left(\frac{1}{2}u^{2}-uv_{h} \right)~dxdt- \int_{Q_1}\xi^{\alpha}\psi_{\varepsilon} u\partial_{t}v_{h}~dxdt\\
        &\leq\beta \int_{Q_1}\xi^{\alpha}\psi'_{\varepsilon}\Gamma^{2}_{+}(u)~dxdt\\
        &\underset{\varepsilon\downarrow0}{\longrightarrow}\beta\int_{K_{\rho_1}\times\{t_1\}}\Gamma^{2}_{+}(u)~dx-
        \beta\int_{K_{\rho_2}\times\{t_2\}}\Gamma^{2}_{+}(u)~dx,
    \end{split}
\end{equation}
where the convergence as $\varepsilon\downarrow 0$ follows as in the proof of Lemma 3.1, owing to the distributional convergence of $\psi'_\varepsilon$ to a difference of Dirac delta measures concentrated at $t_1$ and $t_2$. To estimate the remaining terms we let $h\longrightarrow0$ such that
\begin{equation*}
    \frac{\partial}{\partial x_{i}}\left(\xi^{\alpha}(v_{h}-u)\right)\underset{h\downarrow0}{\longrightarrow}- \frac{\partial}{\partial x_{i}}\left(\xi^{\alpha}\beta\left( \Gamma_{+}^{2}\right)'(u)\right)~\text{in}~L^{p_{i}}(\Omega_T),~i=1,..N.
\end{equation*}
Then, we obtain
\begin{equation}
    \begin{split}
        \underset{h\downarrow0}{\lim}&\int_{Q_1}\psi_{\varepsilon}\left|\frac{\partial u}{\partial x_{i}} \right|^{p_{i}-2}\frac{\partial u}{\partial x_{i}}\frac{\partial }{\partial x_{i}}\left(\xi^{\alpha}(v_{h}-u) \right)\\
        \leq&\alpha\beta\int_{Q_1}\xi^{\alpha-1}\psi_{\varepsilon}\left|\frac{\partial u}{\partial x_{i}} \right|^{p_{i}-1}\left|\frac{\partial \xi}{\partial x_{i}} \right|\left(\Gamma^{2}_{+} \right)'(u)~dxdt
        -\beta\int_{Q_1}\xi^{\alpha}\psi_{\varepsilon}\left|\frac{\partial u}{\partial x_{i}} \right|^{p_{i}}\left(\Gamma^{2}_{+} \right)^{''}(u)~dxdt,
    \end{split}
\end{equation}
where the terms linked to $\left(\Gamma^{2}_{+} \right)^{''}$ are well defined since $\frac{\partial u}{\partial x_{i}}=0$ a.e. in the set $\{(u-k)_{+}=C\in\mathbb{R}\}$. Next, by Young's inequality, we get
\begin{equation}
    \begin{split}
       \alpha\beta&\int_{Q_1}\xi^{\alpha-1}\psi_{\varepsilon}\left|\frac{\partial u}{\partial x_{i}} \right|^{p_{i}-1}\left|\frac{\partial \xi}{\partial x_{i}} \right|\left(\Gamma^{2}_{+} \right)'(u)~dxdt\\
        \leq& \nu\int_{Q_1}\psi_{\varepsilon}\xi^{\alpha}\left|\frac{\partial u}{\partial x_{i}} \right|^{p_{i}}2\Gamma_{+}(\Gamma'_{+})^{2}(u)~dxdt+C(\nu)\int_{Q_1}\psi_{\varepsilon}\left|\frac{\partial \xi}{\partial x_{i}} \right|^{p_{i}}2\Gamma_{+}(\Gamma'_{+})^{2-p_{i}}(u)~dxdt,
        \end{split}
\end{equation}
where used (3.14), and took $\alpha$ as in (3.9). Therefore, by putting (3.27) into (3.26), we arrive at
\begin{equation}
    \begin{split}
        \underset{h\downarrow0}{\lim}&\int_{Q_1}\psi_{\varepsilon}\left|\frac{\partial u}{\partial x_{i}} \right|^{p_{i}-2}\frac{\partial u}{\partial x_{i}}\frac{\partial }{\partial x_{i}}\left(\xi^{\alpha}(v_{h}-u) \right)\\
        \leq&-C\int_{Q_1}\xi^{\alpha}\psi_{\varepsilon}\left|\frac{\partial u}{\partial x_{i}} \right|^{p_{i}}2\left(\Gamma'_{+} \right)^{2}~dxdt
        +C\int_{Q_1}\psi_{\varepsilon}\left|\frac{\partial\xi }{\partial x_{i}}\right|^{p_{i}}2\Gamma_{+}\left(\Gamma'_{+} \right)^{2-p_{i}}(u)~dxdt,
    \end{split}
\end{equation}
where we used the fact that $2\Gamma_{+}\left( \Gamma'_{+}\right)^{2}-\left(\Gamma_{+}^{2} \right)^{''}=-2\left( \Gamma'_{+}\right)^{2}$. Hence, by choosing $\nu>0$ small enough so that the first term on the right-hand side of (3.27) can be absorbed into the left-hand side of the energy inequality, and combining (3.25), and (3.28) into (2.3), we obtain (3.16).

In order to prove (3.17), we take as a comparison function
\begin{equation*}\begin{split}
    v_{h}=&[u]_{h}+\beta\left(\Gamma_{-}^{2}\right)'([u]_h)+\left\|\phi-[\phi]_{h} \right\|_{L^{\infty}(Q_1)}\\
    &\geq \phi+\beta\left(\Gamma_{-}^{2}\right)'([u]_h)\geq\phi,
    \end{split}
\end{equation*}
since $\Gamma_{-}$ and $\Gamma_{-}'$ are nonnegative. We then proceed similarly as in the proof of (3.16) to get the desired result.
\end{proof}
\section{Local boundedness of solutions}
\begin{thm}
Under the assumption that
\begin{equation}
    2 < p_{i} < \bar{p}\left(1+\frac{2}{N}\right) \quad \text{for } i=1,\dots,N, \qquad \bar{p} = \left(\frac{1}{N}\sum_{i=1}^{N}\frac{1}{p_{i}}\right)^{-1} < N,
    \label{eq:assump}
\end{equation}
and that the obstacle $\phi \in C^{0}(\Omega_{T})$, any local weak solution $u \in K_\phi(\Omega_T)$ to the obstacle problem related to the main equation in the sense of Definition 2.3 is locally bounded.
\end{thm}

\begin{proof}
Fix $\rho>0$ sufficiently small such that
\[
Q(\rho^{\bar p},\rho):=K_\rho\times(-\rho^{\bar p},0)
\subset\Omega_T,
\qquad
K_\rho:=\prod_{i=1}^N
\left(-\rho^{\bar p/p_i},\rho^{\bar p/p_i}\right).
\]
Let
\[
\frac12\rho\le \rho_2<\rho_1\le\rho,
\qquad
-\rho^{\bar p}\le -s_1<-s_2\le-\frac12\rho^{\bar p}.
\]
Choose a smooth cutoff function $0\le\xi\le1$ such that
$\xi=1$ in $Q(s_2,\rho_2)$ and $\xi=0$ on the parabolic boundary of
$Q(s_1,\rho_1)$. Then
\begin{equation}
\label{eq:cutoff_bounds}
\left|\frac{\partial\xi}{\partial x_i}\right|
\le
\frac{C}{\rho_1^{\bar p/p_i}-\rho_2^{\bar p/p_i}},
\qquad
\left|\frac{\partial\xi}{\partial t}\right|
\le
\frac{C}{s_1-s_2}.
\end{equation}

Set
\[
l=\bar p\left(1+\frac2N\right)
\]
and choose $\alpha\ge p^+$. For
\[
w=(u-k)_+\xi^\alpha,
\]
we have
\[
\left(
\int_{Q(s_2,\rho_2)\cap\{u>k\}}
(u-k)^l\,dx\,dt
\right)^{\frac{N}{N+\bar p}}
\le
\left(
\int_{Q(s_1,\rho_1)}w^l\,dx\,dt
\right)^{\frac{N}{N+\bar p}}.
\]
Applying Lemma 2.2 to $w$ yields
\begin{equation}
\label{eq:lemma22}
\left(
\int_{Q(s_1,\rho_1)}w^l\,dx\,dt
\right)^{\frac{N}{N+\bar p}}
\le
C\left(
\operatorname*{ess\,sup}_{-s_1<t<0}
\int_{K_{\rho_1}}w^2(x,t)\,dx
+
\sum_{i=1}^N
\int_{Q(s_1,\rho_1)}
|\partial_iw|^{p_i}\,dx\,dt
\right).
\end{equation}

Since
\[
\partial_iw
=
\xi^\alpha\partial_i(u-k)_+
+
\alpha(u-k)_+\xi^{\alpha-1}\partial_i\xi,
\]
the inequality
$|a+b|^{p_i}\le 2^{p_i-1}(|a|^{p_i}+|b|^{p_i})$ and $\xi\le1$ give
\[
|\partial_iw|^{p_i}
\le
C\left(
\xi^{\alpha p_i}|\partial_i u|^{p_i}\mathbf 1_{\{u>k\}}
+
|\partial_i\xi|^{p_i}(u-k)_+^{p_i}
\right).
\]
Consequently,
\[
\begin{aligned}
\sum_{i=1}^N
\int_{Q(s_1,\rho_1)}
|\partial_iw|^{p_i}\,dx\,dt
\le{}&
C\sum_{i=1}^N
\int_{Q(s_1,\rho_1)\cap\{u>k\}}
\xi^{\alpha p_i}|\partial_i u|^{p_i}\,dx\,dt
\\
&+
C\sum_{i=1}^N
\frac{
\displaystyle
\int_{Q(s_1,\rho_1)\cap\{u>k\}}
(u-k)^{p_i}\,dx\,dt
}{
\left(
\rho_1^{\bar p/p_i}-\rho_2^{\bar p/p_i}
\right)^{p_i}
}.
\end{aligned}
\]

Now let
\[
k>\left(\sup_{Q(s_1,\rho_1)}\phi\right)_+.
\]
By Lemma 3.1,
\begin{equation}
\label{eq:lemma31}
\begin{aligned}
\operatorname*{ess\,sup}_{-s_1<t<0}
&\int_{K_{\rho_1}}
(u-k)_+^2\xi^{2\alpha}\,dx
+
\sum_{i=1}^N
\int_{Q(s_1,\rho_1)\cap\{u>k\}}
\xi^{\alpha p_i}|\partial_i u|^{p_i}\,dx\,dt
\\
&\le
C\int_{Q(s_1,\rho_1)\cap\{u>k\}}
\left(
\left|\frac{\partial\xi}{\partial t}\right|(u-k)^2
+
\sum_{i=1}^N
|\partial_i\xi|^{p_i}(u-k)^{p_i}
\right)\,dx\,dt
\\
&\le
C\frac{1}{s_1-s_2}
\int_{Q(s_1,\rho_1)\cap\{u>k\}}
(u-k)^2\,dx\,dt
\\
&\quad+
C\sum_{i=1}^N
\frac{
\displaystyle
\int_{Q(s_1,\rho_1)\cap\{u>k\}}
(u-k)^{p_i}\,dx\,dt
}{
\left(
\rho_1^{\bar p/p_i}-\rho_2^{\bar p/p_i}
\right)^{p_i}
}.
\end{aligned}
\end{equation}

Combining \eqref{eq:lemma22} with \eqref{eq:lemma31}, we obtain
\begin{equation}
\label{eq:sobolev_energy}
\begin{aligned}
\left(
\int_{Q(s_2,\rho_2)\cap\{u>k\}}
(u-k)^l\,dx\,dt
\right)^{\frac{N}{N+\bar p}}
\le
C\sum_{i=1}^N
\bigg\{
&\frac{1}{s_1-s_2}
\int_{Q(s_1,\rho_1)\cap\{u>k\}}
(u-k)^2\,dx\,dt
\\
&+
\frac{
\displaystyle
\int_{Q(s_1,\rho_1)\cap\{u>k\}}
(u-k)^{p_i}\,dx\,dt
}{
\left(
\rho_1^{\bar p/p_i}-\rho_2^{\bar p/p_i}
\right)^{p_i}
}
\bigg\}.
\end{aligned}
\end{equation}

Let
\[
k>h>\left(\sup_{Q(s_1,\rho_1)}\phi\right)_+,
\]
then,
\begin{equation*}
    \begin{split}
\left| Q(s_{1},\rho_{1})\cap\{u>k\} \right|&=\int_{ Q(s_{1},\rho_{1})\cap\{u>k\}}\left|\frac{k-h}{k-h}   \right|^{l}~dxdt
\leq \int_{ Q(s_{1},\rho_{1})\cap\{u>k\}}\left|\frac{u-h}{k-h}   \right|^{l}~dxdt\\
&\leq \int_{ Q(s_{1},\rho_{1})\cap\{u>h\}}\left|\frac{u-h}{k-h}   \right|^{l}~dxdt.
\end{split}
\end{equation*}
Since $l>2$, Hölder's inequality gives
\[
\begin{aligned}
\int_{Q(s_1,\rho_1)\cap\{u>k\}}
(u-k)^2\,dx\,dt
&\le
\left(
\int_{Q(s_1,\rho_1)\cap\{u>h\}}
(u-h)^l\,dx\,dt
\right)^{2/l}
|Q(s_1,\rho_1)\cap\{u>k\}|^{1-2/l}
\\
&\le
\frac{1}{(k-h)^{l-2}}
\int_{Q(s_1,\rho_1)\cap\{u>h\}}
(u-h)^l\,dx\,dt.
\end{aligned}
\]
Likewise, since $p_i<l$,
\[
\int_{Q(s_1,\rho_1)\cap\{u>k\}}
(u-k)^{p_i}\,dx\,dt
\le
\frac{1}{(k-h)^{l-p_i}}
\int_{Q(s_1,\rho_1)\cap\{u>h\}}
(u-h)^l\,dx\,dt.
\]
Hence \eqref{eq:sobolev_energy} gives
\begin{equation}
\label{eq:combined}
\begin{aligned}
\left(
\int_{Q(s_2,\rho_2)\cap\{u>k\}}
(u-k)^l\,dx\,dt
\right)^{\frac{N}{N+\bar p}}
\le&
C
\left(
\int_{Q(s_1,\rho_1)\cap\{u>h\}}
(u-h)^l\,dx\,dt
\right)\\
&\times
\sum_{i=1}^N
\bigg\{
\frac{1}{(s_1-s_2)(k-h)^{l-2}}
\\
&+
\frac{1}{
\left(
\rho_1^{\bar p/p_i}-\rho_2^{\bar p/p_i}
\right)^{p_i}
(k-h)^{l-p_i}}
\bigg\}.
\end{aligned}
\end{equation}

By the absolute continuity of the Lebesgue integral, we can choose
\[
H>\left(\sup_{Q(\rho^{\bar p},\rho)}\phi\right)_+
\]
so large that
\begin{equation}
\label{eq:abs_cont}
\int_{-\rho^{\bar p}}^0
\int_{K_\rho}
(u-H)_+^l\,dx\,dt
\le
\varepsilon\rho^{N+\bar p},
\end{equation}
where $\varepsilon>0$ will be chosen below.

For $m=0,1,2,\dots$, define
\[
k_m
=
2H-\frac{H}{2^m},
\qquad
\rho_m
=
\left(\frac12+\frac{1}{2^{m+1}}\right)\rho,
\]
and
\[
s_m
=
\left(\frac12+\frac{1}{2^{m+1}}\right)\rho^{\bar p},
\qquad
Q_m
=
K_{\rho_m}\times(-s_m,0).
\]
Set
\[
J_m
=
\int_{Q_m\cap\{u>k_m\}}
(u-k_m)^l\,dx\,dt.
\]

We apply \eqref{eq:combined} with
\[
Q(s_2,\rho_2)=Q_{m+1},
\qquad
Q(s_1,\rho_1)=Q_m,
\qquad
k=k_{m+1},
\qquad
h=k_m.
\]
Observe that
\[
s_m-s_{m+1}
=
\frac{\rho^{\bar p}}{2^{m+2}},
\qquad
k_{m+1}-k_m
=
\frac{H}{2^{m+1}}.
\]
Moreover, since
\[
\rho_{m+1}\ge\frac{\rho}{2},
\qquad
\rho_m-\rho_{m+1}=\frac{\rho}{2^{m+2}},
\]
the mean value theorem gives
\[
\rho_m^{\bar p/p_i}
-
\rho_{m+1}^{\bar p/p_i}
\ge
c\,2^{-m}\rho^{\bar p/p_i}
\]
with a constant $c>0$ independent of $m$. Therefore,
\[
\frac{1}{(s_m-s_{m+1})(k_{m+1}-k_m)^{l-2}}
\le
\frac{C\,2^{m(l-1)}}{\rho^{\bar p}H^{l-2}},
\]
and
\[
\frac{1}{
\left(
\rho_m^{\bar p/p_i}-\rho_{m+1}^{\bar p/p_i}
\right)^{p_i}
(k_{m+1}-k_m)^{l-p_i}}
\le
\frac{C\,2^{ml}}{\rho^{\bar p}H^{l-p_i}}.
\]
Thus, with $d=2^{l+1}$ and after enlarging $C$ if necessary,
\begin{equation}
\label{eq:recurrence_full}
J_{m+1}^{\frac{N}{N+\bar p}}
\le
C\left(
\frac{d^m}{\rho^{\bar p}H^{l-2}}
+
\sum_{i=1}^N
\frac{d^m}{\rho^{\bar p}H^{l-p_i}}
\right)J_m.
\end{equation}

Since $H>1$ and $p_i<l$, we have
\[
H^{-(l-2)}\le1,
\qquad
H^{-(l-p_i)}\le1.
\]
Hence, after absorbing the finite sum into $C$,
\[
J_{m+1}^{\frac{N}{N+\bar p}}
\le
\frac{C d^m}{\rho^{\bar p}}J_m.
\]
Writing
\[
J_m
=
J_m^{\frac{N}{N+\bar p}}
J_m^{\frac{\bar p}{N+\bar p}},
\]
we obtain
\begin{equation}
\label{eq:recurrence_simplified}
J_{m+1}^{\frac{N}{N+\bar p}}
\le
C J_m^{\frac{N}{N+\bar p}}
\left(
\frac{d^m}{\rho^{\bar p}}
J_m^{\frac{\bar p}{N+\bar p}}
\right).
\end{equation}

From \eqref{eq:abs_cont},
\[
J_0
=
\int_{Q_0\cap\{u>H\}}
(u-H)^l\,dx\,dt
\le
\varepsilon\rho^{N+\bar p}.
\]
We claim that there exists $\delta\in(0,1)$ such that
\begin{equation}
\label{eq:induction_hyp}
J_m
\le
\delta^m\varepsilon\rho^{N+\bar p},
\qquad
m=0,1,2,\dots.
\end{equation}
The case $m=0$ follows from the preceding estimate. Assume that \eqref{eq:induction_hyp} holds for some $m\ge0$. By \eqref{eq:recurrence_simplified},
\[
\begin{aligned}
J_{m+1}^{\frac{N}{N+\bar p}}
&\le
C J_m^{\frac{N}{N+\bar p}}
\left(
\frac{d^m}{\rho^{\bar p}}
\left(
\delta^m\varepsilon\rho^{N+\bar p}
\right)^{\frac{\bar p}{N+\bar p}}
\right)
\\
&=
J_m^{\frac{N}{N+\bar p}}
\left(
C d^m
\delta^{\frac{m\bar p}{N+\bar p}}
\varepsilon^{\frac{\bar p}{N+\bar p}}
\right).
\end{aligned}
\]
Choose $\delta\in(0,1)$ such that
\[
d\,\delta^{\frac{\bar p}{N+\bar p}}\le1.
\]
Then choose $\varepsilon>0$ sufficiently small that
\begin{equation}
\label{eq:conditions}
C\varepsilon^{\frac{\bar p}{N+\bar p}}
\le
\delta^{\frac{N}{N+\bar p}}.
\end{equation}
It follows that
\[
C d^m
\delta^{\frac{m\bar p}{N+\bar p}}
\varepsilon^{\frac{\bar p}{N+\bar p}}
\le
\delta^{\frac{N}{N+\bar p}},
\]
and consequently
\[
\begin{aligned}
J_{m+1}^{\frac{N}{N+\bar p}}
&\le
J_m^{\frac{N}{N+\bar p}}
\delta^{\frac{N}{N+\bar p}}
\\
&\le
\left(
\delta^m\varepsilon\rho^{N+\bar p}
\right)^{\frac{N}{N+\bar p}}
\delta^{\frac{N}{N+\bar p}}
\\
&=
\left(
\delta^{m+1}\varepsilon\rho^{N+\bar p}
\right)^{\frac{N}{N+\bar p}}.
\end{aligned}
\]
Thus
\[
J_{m+1}
\le
\delta^{m+1}\varepsilon\rho^{N+\bar p},
\]
and the induction is complete.

Since $0<\delta<1$, we have $J_m\to0$. Moreover,
\[
k_m\to2H,
\qquad
\rho_m\to\frac{\rho}{2},
\qquad
s_m\to\frac{\rho^{\bar p}}2.
\]
For almost every $(x,t)$,
\[
(u-k_m)_+^l\mathbf 1_{Q_m}
\longrightarrow
(u-2H)_+^l
\mathbf 1_{K_{\rho/2}\times(-\rho^{\bar p}/2,0)}.
\]
Furthermore,
\[
0\le
(u-k_m)_+^l\mathbf1_{Q_m}
\le
(u-H)_+^l\mathbf1_{Q_0},
\]
and the right-hand side is integrable by \eqref{eq:abs_cont}. Hence the dominated convergence theorem gives
\[
\begin{aligned}
0
=
\lim_{m\to\infty}J_m
&=
\int_{K_{\rho/2}\times(-\rho^{\bar p}/2,0)}
(u-2H)_+^l\,dx\,dt.
\end{aligned}
\]
Therefore,
\[
\operatorname*{ess\,sup}_{K_{\rho/2}\times(-\rho^{\bar p}/2,0)}
u
\le
2H.
\]
This with $u\geq\phi$ gives the local boundedness of $u$ over $\Omega_{T}$.
\end{proof}

\section{Toward H\"older continuity}
Let $R\in(0,1)$ be small enough such that the standard box cylinder
\[
Q(R^{2},R)
=
(-R,R)^{N}\times(-R^{2},0)
\subset\Omega_{T},
\]
where $(-R,R)^{N}$ is the $N$-dimensional spatial box of half-length $R$, and where $u$ is locally bounded by Theorem 4.1. We assume further that $\phi\in C^{0;\beta,\frac{\beta}{2}}(\Omega_{T})$ for some
$\beta\in(0,1)$, with
\begin{equation*}
    [\phi]_{C^{0;\beta,\frac{\beta}{2}}}
    :=\sup_{(x,t),(y,s)\in\Omega_{T}}
    \frac{|\phi(x,t)-\phi(y,s)|}{|x-y|^{\beta}+|t-s|^{\frac{\beta}{2}}}.
\end{equation*}
For some $\lambda>1$ and $\vartheta\in(0,\beta)$ to be fixed later, define
\begin{equation}
    H(\rho):=\max\left\{2^{\lambda}\rho^{\vartheta},\,
    2\,\operatorname*{ess\,osc}_{Q(\rho^{2},\rho)}\phi\right\},
    \qquad 0\le\rho\le R. \label{eq:H}
\end{equation}
Then $H$ is continuous and increasing. Since $\phi$ is H\"older continuous,
\begin{equation}
    \operatorname*{ess\,osc}_{Q(\rho^{2},\rho)}\phi
    \le [\phi]_{C^{0;\beta,\frac{\beta}{2}}}\rho^{\beta},
    \qquad 0\le\rho\le R. \label{eq:osc-phi}
\end{equation}
If
\[
\operatorname*{ess\,osc}_{Q(\rho^{2},\rho)}u\le H(\rho)
\quad\text{for all }\rho\in(0,R],
\]
then $u$ is H\"older continuous. Assume this fails. Then there exists
$\rho_{0}\in(0,R]$ such that
\begin{equation}
    H(\rho_{0})\le \operatorname*{ess\,osc}_{Q(\rho_{0}^{2},\rho_{0})}u,
    \qquad
    \operatorname*{ess\,osc}_{Q(\rho^{2},\rho)}u
    \le C R^{-\beta}H(\rho)
    \quad \forall\,\rho\in[\rho_{0},R]. \label{eq:rho0}
\end{equation}
Set
\[
\mu^{+}=\operatorname*{ess\,sup}_{Q(\rho_{0}^{2},\rho_{0})}u,\quad
\mu^{-}=\operatorname*{ess\,inf}_{Q(\rho_{0}^{2},\rho_{0})}u,\quad
\omega=\mu^{+}-\mu^{-}.
\]
Denote
\[
p^{+}=\max\{p_{1},\dots,p_{N}\},\qquad
p^{-}=\min\{p_{1},\dots,p_{N}\},\qquad
\frac{1}{\bar p}=\frac{1}{N}\sum_{i=1}^{N}\frac{1}{p_{i}}.
\]

For $\rho>0$ and $\omega>0$, define the anisotropic spatial box
\[
K_{\rho}
:=\prod_{i=1}^{N}
\left(
-\rho^{\frac{p^{+}}{p_{i}}}\omega^{\frac{p_{i}-p^{+}}{p_{i}}},
 \rho^{\frac{p^{+}}{p_{i}}}\omega^{\frac{p_{i}-p^{+}}{p_{i}}}
\right),
\]
and the intrinsic cylinder
\begin{equation}
    \mathcal{Q}_{a,\omega}(\rho)
    := K_{\rho}\times(-a\rho^{p^{+}},0),
    \qquad a>0.
\end{equation}
Its time length is $a\rho^{p^{+}}$, and its spatial half-length in the
$x_{i}$-direction is
\[
L_{i}(\rho)=\rho^{\frac{p^{+}}{p_{i}}}\omega^{\frac{p_{i}-p^{+}}{p_{i}}}.
\]

Set
\begin{equation}
    a_{0}=\left(\frac{\omega}{2^{\lambda}}\right)^{2-p^{+}},
    \qquad
    \widehat Q_{0}:=\mathcal{Q}_{a_{0},\omega}(\rho_{0}).
    \label{eq:a0}
\end{equation}
We prove that
\begin{equation}
    \widehat Q_{0}\subset Q(\rho_{0}^{2},\rho_{0}). \label{eq:Q0-inclusion}
\end{equation}
Since $H(\rho_{0})\le\omega$ by \eqref{eq:rho0}, and
$H(\rho_{0})\ge 2^{\lambda}\rho_{0}^{\vartheta}$ by \eqref{eq:H}, we have
\[
2^{\lambda}\rho_{0}^{\vartheta}\le\omega.
\]
Because $\rho_{0}<1$ and $\vartheta<1$, it follows that
$\rho_{0}<\rho_{0}^{\vartheta}$. Hence
\[
2^{\lambda}\rho_{0}<\omega.
\]
Equivalently,
\[
\frac{\omega}{2^{\lambda}}>\rho_{0}.
\]
Since $p^{+}\ge2$, the exponent $2-p^{+}$ is non-positive. Raising the previous
inequality to the power $2-p^{+}$ reverses the inequality:
\[
a_{0}=\left(\frac{\omega}{2^{\lambda}}\right)^{2-p^{+}}
\le \rho_{0}^{2-p^{+}}.
\]
Therefore
\[
a_{0}\rho_{0}^{p^{+}}\le \rho_{0}^{2},
\]
so the time interval of $\widehat Q_{0}$ is contained in $(-\rho_{0}^{2},0)$.

For the spatial inclusion, for each $i$,
\[
\rho_{0}^{\frac{p^{+}}{p_{i}}}\omega^{\frac{p_{i}-p^{+}}{p_{i}}}\le\rho_{0}
\]
is equivalent to
\[
\rho_{0}^{\frac{p^{+}}{p_{i}}-1}
\le \omega^{\frac{p^{+}-p_{i}}{p_{i}}}.
\]
If $p_{i}=p^{+}$, this is trivial. If $p_{i}<p^{+}$, the exponent
$\frac{p^{+}-p_{i}}{p_{i}}$ is positive, and the inequality is equivalent to
$\rho_{0}\le\omega$. This is true because
\[
\omega>2^{\lambda}\rho_{0}^{\vartheta}>2^{\lambda}\rho_{0}>\rho_{0}.
\]
Hence $K_{\rho_{0}}\subset B_{\rho_{0}}$, and \eqref{eq:Q0-inclusion} follows.
Therefore
\begin{equation}
\operatorname*{ess\,osc}_{\widehat Q_{0}}u\le\omega.
\end{equation}
Moreover,
\begin{equation}
\begin{aligned}
\operatorname*{ess\,sup}_{\widehat Q_{0}}\phi
&\le \operatorname*{ess\,sup}_{Q(\rho_{0}^{2},\rho_{0})}\phi \\
&= \operatorname*{ess\,inf}_{Q(\rho_{0}^{2},\rho_{0})}\phi
   +\operatorname*{ess\,osc}_{Q(\rho_{0}^{2},\rho_{0})}\phi \\
&\le \mu^{-}+\frac{\omega}{2}.
\end{aligned}
\end{equation}
Thus for every $k\ge \mu^{-}+\omega/2$ we have
$k\ge\operatorname*{ess\,sup}_{\widehat Q_{0}}\phi$.

Next, let
\[
    \theta=\left(\frac{\omega}{2}\right)^{2-p^{+}}.
\]
Then
\[
\frac{a_{0}}{\theta}
=
\left(\frac{\omega}{2^{\lambda}}\right)^{2-p^{+}}
\left(\frac{\omega}{2}\right)^{p^{+}-2}
=
2^{(\lambda-1)(p^{+}-2)}.
\]
Since $\lambda>1$, we have $a_{0}>\theta$. Choose $\lambda$ so that
$N_{0}:=2^{(\lambda-1)(p^{+}-2)}$ is an integer.

For any $\tau^{*}$ satisfying
\begin{equation}
    -(a_{0}-\theta)\rho_{0}^{p^{+}}<\tau^{*}<0, \label{eq:tau-range}
\end{equation}
consider the subcylinder
\[
Q^{*}=(0,\tau^{*})+\mathcal{Q}_{\theta,\omega}(\rho_{0}).
\]
Its time interval is
\[
(\tau^{*}-\theta\rho_{0}^{p^{+}},\,\tau^{*}).
\]
Because of \eqref{eq:tau-range}, its bottom satisfies
\[
\tau^{*}-\theta\rho_{0}^{p^{+}}>-a_{0}\rho_{0}^{p^{+}},
\]
and its top satisfies $\tau^{*}<0$. Hence
\[
Q^{*}\subset \widehat Q_{0}.
\]
Moreover, the cylinders
\[
(0,-k\theta\rho_{0}^{p^{+}})+\mathcal{Q}_{\theta,\omega}(\rho_{0}),
\qquad k=0,1,\dots,N_{0}-1,
\]
tile $\widehat Q_{0}$ exactly in the time direction.

Now fix $\nu_{0}\in(0,1)$ to be chosen later depending only on the data.
We distinguish two alternatives. First, there exists a cylinder
$Q^{*}=(0,\tau^{*})+\mathcal{Q}_{\theta,\omega}(\rho_{0})$ as above such that
\begin{equation}
    \left|\left\{(x,t)\in Q^{*}:
    u<\mu^{-}+\frac{\omega}{2}\right\}\right|
    \le \nu_{0}\,|Q^{*}|. \label{eq:first-alt}
\end{equation}
Second, if \eqref{eq:first-alt} fails, then for every such cylinder
$Q^{*}=(0,\tau^{*})+\mathcal{Q}_{\theta,\omega}(\rho_{0})$ we have
\begin{equation}
    \left|\left\{(x,t)\in Q^{*}:
    u\ge\mu^{-}+\frac{\omega}{2}\right\}\right|
    \le (1-\nu_{0})\,|Q^{*}|. \label{eq:second-alt}
\end{equation}

\subsection{First alternative}
Assume that \eqref{eq:first-alt} holds.

\begin{lem}
    There exists a constant $\nu_{0}\in(0,1)$, depending only on
    $N,\{p_{i}\}_{i=1}^{N},\bar p$, such that
    \begin{equation}
        u>\mu^{-}+\frac{\omega}{4}
        \quad\text{a.e. in}\quad
        Q_{\infty}:=(0,\tau^{*})+\mathcal{Q}_{\theta,\omega}\!\left(\frac{\rho_{0}}{2}\right).
        \label{eq:first-conclusion}
    \end{equation}
\end{lem}

\begin{proof}
Define
\[
\rho_n=\frac{\rho_0}{2}+\frac{\rho_0}{2^{n+1}},
\qquad
k_n=\mu^-+\frac{\omega}{4}+\frac{\omega}{2^{n+2}},
\qquad
Q_n=(0,\tau^*)+\mathcal{Q}_{\theta,\omega}(\rho_n),
\]
and
\[
A_n=Q_n\cap\{u<k_n\}.
\]
Let
\[
w_n=(u-k_n)_- .
\]
Choose a cutoff function $\zeta_n$ such that
\[
\zeta_n=1\quad\text{on }Q_{n+1},
\qquad
\zeta_n=0\quad\text{on the parabolic boundary of }Q_n,
\]
and
\[
|\partial_{x_i}\zeta_n|
\le C2^n\rho_0^{-\frac{p^{+}}{p_i}}\omega^{\frac{p^{+}-p_i}{p_i}},
\qquad
|\partial_t\zeta_n|
\le C2^n\rho_0^{-p^{+}}\omega^{p^{+}-2}.
\]

Apply the energy estimate (3.1) to $w_n$ with test function
$\zeta_n^\alpha$, where $\alpha\ge p^{+}$. Since $\zeta_n=0$ at the
bottom of $Q_n$, the initial term vanishes. Hence
\[
\begin{aligned}
&\sup_{t\in(\tau^*-\theta\rho_n^{p^{+}},\tau^*)}
\int_{K_{\rho_n}}\zeta_n^\alpha w_n^2\,dx
+\sum_{i=1}^N\iint_{Q_n}\zeta_n^\alpha|\partial_{x_i}w_n|^{p_i}\,dxdt\\
&\le C\sum_{i=1}^N
\biggl\{
\iint_{Q_n} w_n^2|\partial_t\zeta_n|\,dxdt
+\iint_{Q_n} w_n^{p_i}|\partial_{x_i}\zeta_n|^{p_i}\,dxdt
\biggr\}.
\end{aligned}
\]
On $A_n$,
\[
w_n\le k_n-\mu^-\le \frac{\omega}{2}.
\]
Thus $w_n^2\le \omega^2/4$ and $w_n^{p_i}\le \omega^{p_i}/2^{p_i}$.
Using the derivative bounds for $\zeta_n$, we get
\[
\iint_{Q_n} w_n^2|\partial_t\zeta_n|\,dxdt
\le C2^n\rho_0^{-p^{+}}\omega^{p^{+}}|A_n|,
\]
and, for each $i$,
\[
\iint_{Q_n} w_n^{p_i}|\partial_{x_i}\zeta_n|^{p_i}\,dxdt
\le C2^{np_i}\rho_0^{-p^{+}}\omega^{p^{+}}|A_n|
\le C2^{np^{+}}\rho_0^{-p^{+}}\omega^{p^{+}}|A_n|.
\]
Therefore
\begin{equation}
\begin{aligned}
&\sup_{t\in(\tau^*-\theta\rho_n^{p^{+}},\tau^*)}
\int_{K_{\rho_n}}\zeta_n^\alpha w_n^2\,dx
+\sum_{i=1}^N\iint_{Q_n}\zeta_n^\alpha|\partial_{x_i}w_n|^{p_i}\,dxdt\\
&\le C2^{np^{+}}\rho_0^{-p^{+}}\omega^{p^{+}}|A_n|.
\end{aligned}
\label{eq:weighted-energy}
\end{equation}

Now define
\[
v_n:=w_n\zeta_n^{\alpha/2}.
\]
Since $\zeta_n=0$ on the parabolic boundary of $Q_n$, the function
$v_n$ vanishes there. Moreover,
\[
v_n^2=w_n^2\zeta_n^\alpha,
\]
so
\[
\sup_{t\in(\tau^*-\theta\rho_n^{p^{+}},\tau^*)}
\int_{K_{\rho_n}}v_n^2\,dx
=
\sup_{t\in(\tau^*-\theta\rho_n^{p^{+}},\tau^*)}
\int_{K_{\rho_n}}\zeta_n^\alpha w_n^2\,dx.
\]
For the gradient, by the product rule,
\[
\partial_{x_i}v_n
=
\zeta_n^{\alpha/2}\partial_{x_i}w_n
+
\frac\alpha2 w_n\zeta_n^{\alpha/2-1}\partial_{x_i}\zeta_n.
\]
Hence
\[
|\partial_{x_i}v_n|^{p_i}
\le
C\zeta_n^{\alpha p_i/2}|\partial_{x_i}w_n|^{p_i}
+
C w_n^{p_i}\zeta_n^{(\alpha/2-1)p_i}|\partial_{x_i}\zeta_n|^{p_i}.
\]
Since $p_i\ge2$, we have $\alpha p_i/2\ge\alpha$, and because
$0\le\zeta_n\le1$,
\[
\zeta_n^{\alpha p_i/2}\le \zeta_n^\alpha.
\]
Thus the first term is controlled by the weighted gradient in
\eqref{eq:weighted-energy}. For the second term, using $w_n\le\omega/2$
and the bound for $\partial_{x_i}\zeta_n$,
\[
\iint_{Q_n} w_n^{p_i}\zeta_n^{(\alpha/2-1)p_i}|\partial_{x_i}\zeta_n|^{p_i}\,dxdt
\le C2^{np^{+}}\rho_0^{-p^{+}}\omega^{p^{+}}|A_n|.
\]
Combining these estimates, we obtain
\[
\sup_{t\in(\tau^*-\theta\rho_n^{p^{+}},\tau^*)}
\int_{K_{\rho_n}}v_n^2\,dx
+\sum_{i=1}^N\iint_{Q_n}|\partial_{x_i}v_n|^{p_i}\,dxdt
\le C2^{np^{+}}\rho_0^{-p^{+}}\omega^{p^{+}}|A_n|.
\]

By Lemma 2.2 with $l=\bar p(1+2/N)$,
\[
\iint_{Q_n}|v_n|^l\,dxdt
\le C\left(2^{np^{+}}\rho_0^{-p^{+}}\omega^{p^{+}}|A_n|\right)^{\frac{N+\bar p}{N}}.
\]

Since $\zeta_n=1$ on $Q_{n+1}$, we have
$v_n=w_n=(u-k_n)_-$ a.e. in $Q_{n+1}$.
On $A_{n+1}$,
\[
(u-k_n)_-\ge k_n-k_{n+1}
=\frac{\omega}{2^{n+3}}.
\]
Therefore
\[
|A_{n+1}|\left(\frac{\omega}{2^{n+3}}\right)^l
\le C\left(2^{np^{+}}\rho_0^{-p^{+}}\omega^{p^{+}}|A_n|\right)^{\frac{N+\bar p}{N}}.
\]
Rearranging,
\[
|A_{n+1}|
\le C2^{n(l+p^{+}\frac{N+\bar p}{N})}
\rho_0^{-p^{+}\frac{N+\bar p}{N}}
\omega^{p^{+}\frac{N+\bar p}{N}-l}
|A_n|^{\frac{N+\bar p}{N}}.
\]

Now define
\[
X_n=\frac{|A_n|}{|Q_n|}.
\]
The volume of $Q_n$ is
\[
|Q_n|=C\rho_n^A\omega^B,
\qquad
A=p^{+}\Bigl(1+\frac{N}{\bar p}\Bigr),\quad
B=2+N-p^{+}\Bigl(1+\frac{N}{\bar p}\Bigr).
\]
Since $\rho_n\sim\rho_0$, dividing the previous inequality by $|Q_{n+1}|$
and using $|A_n|=X_n|Q_n|$ gives
\[
X_{n+1}
\le C 2^{n(l+p^{+}\frac{N+\bar p}{N})}
X_n^{1+\frac{\bar p}{N}},
\]
where $C$ depends only on the data. The powers of $\rho_0$ and
$\omega$ cancel exactly.

Let
\[
\alpha_0=\frac{\bar p}{N},
\qquad
B_0=2^{l+p^{+}\frac{N+\bar p}{N}}.
\]
Then
\[
X_{n+1}\le C~B_0^nX_n^{1+\alpha_0}.
\]
By Lemma 2.5, $X_n\to0$ if
\[
X_0\le \nu_0:=C^{-N/\bar p}B_0^{-N^2/\bar p^2}.
\]
The first alternative \eqref{eq:first-alt} gives exactly this bound on
$X_0$. Since $k_n\downarrow \mu^-+\omega/4$, we conclude
\[
u>\mu^-+\frac{\omega}{4}
\quad\text{a.e. in }Q_\infty.
\]
\end{proof}
Now, our next goal is to have similar estimations in smaller cylinders. Consequently, let
\begin{equation}
    s:=\theta\left(\frac{\rho_{0}}{2}\right)^{p^{+}},
    \qquad
    -t_{0}:=\tau^{*}-s.
\end{equation}
Then
\[
(0,\tau^{*})+\mathcal{Q}_{\theta,\omega}\!\left(\frac{\rho_{0}}{2}\right)
=
K_{\frac{\rho_{0}}{2}}\times(-t_{0},\tau^{*}),
\]
and Lemma 5.1 gives
\begin{equation}
    u(\cdot,-t_{0})>\mu^{-}+\frac{\omega}{4}
    \quad\text{a.e. in }K_{\frac{\rho_{0}}{2}}.
\end{equation}

\begin{lem}
    For every $\tilde{\nu}\in(0,1)$ there exists $n_{1}\in\mathbb N$,
    depending only on the data and $\tilde{\nu}$, such that
    \begin{equation}
        \left|\left\{x\in K_{\frac{\rho_{0}}{4}}:
        u(x,t)<\mu^{-}+\frac{\omega}{2^{n_{1}}}\right\}\right|
        \le \tilde{\nu}\,\left|K_{\frac{\rho_{0}}{4}}\right|,
        \qquad \forall t\in(-t_{0},0).
        \label{eq:log-slice}
    \end{equation}
\end{lem}
\begin{proof}
We apply the logarithmic estimate (3.17) to the cylinder $K_{\frac{\rho_{0}}{2}}\times(-t_{0},0)$ and the smaller box $K_{\frac{\rho_{0}}{4}}$. 

For the chosen level $k=\mu^{-}+\frac{\omega}{4}$ and the constant $c=\frac{\omega}{2^{n+2}}$, we have
\[
    k-u \leq H^{-}_{k} = \operatorname*{ess\,sup}_{Q_{\infty}} \left|\left(u-\mu^{-}-\frac{\omega}{4}\right)_{-}\right| \leq \frac{\omega}{4}.
\]

Assuming $H^{-}_{k}>\frac{\omega}{8}$ (otherwise the result is trivial), we can bound the logarithmic function and its derivative. By definition, we have:
\[
    \Gamma_{-} \leq n\ln 2 \quad\text{since}\quad \frac{H^{-}_{k}}{H^{-}_{k}-(u-k)_-+c} \leq \frac{\frac{\omega}{4}}{c} = 2^{n},
\]
\[
    0 \leq \Gamma'_{-} \leq \frac{1}{c} \quad\text{for } u \neq -k+c.
\]

Next, we bound the derivative $|\Gamma'_{-}| = \frac{1}{H^{-}_{k}-(u-k)_-+c}$. The denominator is bounded from above by:
\[
    H^{-}_{k} + u-k- + c \leq \frac{\omega}{4} + 0 + \frac{\omega}{4} = \frac{\omega}{2}.
\]
This implies the derivative is bounded from below:
\[
    |\Gamma'_{-}| \geq \frac{2}{\omega}.
\]
Because $p_{i} \geq 2$, the exponent $2-p_{i}$ is non-positive. Raising both sides to this power reverses the inequality, yielding:
\[
    |\Gamma'_{-}|^{2-p_{i}} \leq \left(\frac{2}{\omega}\right)^{2-p_{i}} = \left(\frac{\omega}{2}\right)^{p_{i}-2}, \quad \text{for } i=1,\dots,N.
\]

In lemma 3.2, the geometric factor $\frac{1}{(\rho_1 - \rho_2)^{p_i}}$ represents the inverse distance between the domain boundaries to the power $p_i$. For our nested anisotropic boxes $K_{\frac{\rho_0}{4}} \subset K_{\frac{\rho_0}{2}}$, the distance in the $i$-th direction is the difference in their half-lengths:
\[
    L_i\left(\frac{\rho_0}{2}\right) - L_i\left(\frac{\rho_0}{4}\right) 
    = \left[ \left(\frac{1}{2}\right)^{\frac{p^{+}}{p_i}} - \left(\frac{1}{4}\right)^{\frac{p^{+}}{p_i}} \right] \rho_0^{\frac{p^{+}}{p_i}}\omega^{\frac{p_i-p^{+}}{p_i}}.
\]

Applying the logarithmic lemma with this anisotropic spatial factor over the time interval $(t_0, \tau^*)$ gives:
\[
    \operatorname*{ess\,sup}_{t\in(-t_{0},0)} \int_{K_{\frac{\rho_{0}}{4}}}\Gamma_{-}^{2}(u)\,dx \leq \int_{K_{\frac{\rho_{0}}{2}}}\Gamma_{-}^{2}(u(\cdot, -t_0))\,dx + C\sum_{i=1}^{N} \frac{1}{\rho_0^{p^{+}}\omega^{p_i-p^{+}}} \iint_{K_{\frac{\rho_{0}}{2}}\times(-t_{0},0)} \Gamma_{-} (\Gamma'_{-})^{2-p_i} \,dxdt.
\]

By Lemma 5.1, $u(\cdot,-t_{0}) > \mu^{-}+\frac{\omega}{4} = k$ almost everywhere in $K_{\frac{\rho_{0}}{2}}$. This implies $(u-k)_- = 0$ at the initial time $t=-t_0$, so the first integral vanishes:
\[
    [\Gamma_{-}(u)](\cdot,-t_{0}) = 0 \quad \text{a.e. in } K_{\frac{\rho_{0}}{2}}.
\]

For the second integral, we estimate the integrand using our bounds and integrate over the time interval of length $ t_0 \leq a_0\rho_{0}^{p^{+}}=\left(\frac{\omega}{2^{\lambda}}\right)^{2-p^{+}}\rho_{0}^{p^{+}}$. Putting everything together, the spatial and temporal weights cancel exactly:
\begin{equation}
\begin{aligned}
    \operatorname*{ess\,sup}_{t\in(-t_{0},0)} \int_{K_{\frac{\rho_{0}}{4}}}\Gamma_{-}^{2}(u)\,dx 
    &\leq C \sum_{i=1}^{N} \left( \rho_0^{-p^{+}}\omega^{p^{+}-p_i} \right) (n 2^{\lambda(p^{+}-2)}) \left(\frac{\omega}{2}\right)^{p_i-2} |K_{\frac{\rho_{0}}{2}}| \left( \omega^{2-p^{+}} \rho_0^{p^{+}} \right) \\
    &\leq C n 2^{\lambda(p^{+}-2)} \sum_{i=1}^{N} \omega^{(p^{+}-p_i) + (p_i-2) + (2-p^{+})} \rho_0^{-p^{+} + p^{+}} |K_{\frac{\rho_{0}}{2}}| \\
    &\leq C n 2^{\lambda(p^{+}-2)} |K_{\frac{\rho_{0}}{4}}|,
\end{aligned}
\end{equation}
where we used $|K_{\frac{\rho_{0}}{2}}| = 2^{\sum \frac{p^{+}}{p_i}} |K_{\frac{\rho_{0}}{4}}| \leq C |K_{\frac{\rho_{0}}{4}}|$.

To obtain a lower bound on the left-hand side, fix $t\in(-t_{0},0)$ and define the set:
\[
    T_{t} = \left\{ x\in K_{\frac{\rho_{0}}{4}}: u(x,t) < \mu^{-}+\frac{\omega}{2^{n+2}} \right\} \subset K_{\frac{\rho_{0}}{2}}.
\]

On $T_{t}$ we have:
\[
    [\Gamma_{-}(u)]^{2} \geq [\ln 2^{n-1}]^{2} = (n-1)^{2}(\ln 2)^{2},
\]
because:
\[
    \frac{H^{-}_{k}}{H^{-}_{k}+u-k+\frac{\omega}{2^{n+2}}} \geq \frac{\frac{\omega}{4}}{\frac{\omega}{4}+u-k+\frac{\omega}{2^{n+2}}} \geq \frac{\frac{\omega}{4}}{\frac{\omega}{2^{n+1}}} = 2^{n-1}.
\]

Integrating over $T_t$ instead of the full box $K_{\frac{\rho_{0}}{4}}$ gives the lower bound:
\[
    (n-1)^{2}(\ln 2)^{2}|T_{t}| \leq \int_{T_{t}}\Gamma_{-}^{2}(u(x,t))\,dx \leq \int_{K_{\frac{\rho_{0}}{4}}}\Gamma_{-}^{2}(u(x,t))\,dx \leq C n 2^{\lambda(p^{+}-2)} |K_{\frac{\rho_{0}}{4}}|.
\]

Dividing through by $(n-1)^2 (\ln 2)^2$ yields:
\[
    |T_{t}| \leq C\frac{n}{(n-1)^{2}} 2^{\lambda(p^{+}-2)} |K_{\frac{\rho_{0}}{4}}|, \qquad \forall t\in(-t_{0},0).
\]

By choosing $n$ sufficiently large such that
\[
    n > 1 + 2C\frac{2^{\lambda(p^{+}-2)}}{\tilde{\nu}}
\]
and setting $n_{1}=n+2$, we achieve:
\[
    |T_{t}| \leq \tilde{\nu}\,|K_{\frac{\rho_{0}}{4}}|, \qquad \forall t\in(-t_{0},0).
\]

Since $T_{t} = \{x\in K_{\rho_0/4}:u(x,t) < \mu^{-}+\omega/2^{n_{1}}\}$, the claim follows.
\end{proof}

\begin{lem}
There exists \(n_2 \in \mathbb N\), depending only on the data, such that
\begin{equation}
    u > \mu^- + \frac{\omega}{2^{n_2+1}}
    \quad\text{a.e. in}\quad
    K_{\frac{\rho_0}{8}} \times (-t_0, 0).
    \label{eq:first-alt-conclusion}
\end{equation}
\end{lem}

\begin{proof}
Let \(n_2 = n_1\) be as in Lemma 5.2. Define the decreasing radii
\[
\rho_n := \frac{\rho_0}{8} + \frac{\rho_0}{2^{n+3}}, \qquad n = 0,1,2,\dots
\]
and the increasing levels
\[
k_n := \mu^- + \frac{\omega}{2^{n_2+1}} + \frac{\omega}{2^{n_2+1+n}}.
\]
Set
\[
Q_n := K_{\rho_n} \times (-t_0, 0), \qquad A_n := Q_n \cap \{u < k_n\}.
\]
At the initial time \(-t_0\), since
\[
k_0 = \mu^- + \frac{\omega}{2^{n_2}}
\]
and \(u(\cdot, -t_0) > \mu^- + \omega/4\) (by Lemma 5.1), we have
\[
(u - k_n)_-(\cdot, -t_0) = 0 \quad \text{for all } n \ge 0,
\]
provided \(n_2 \ge 2\).

Choose a spatial cutoff \(\zeta_n\) such that
\[
\zeta_n = 1 \quad \text{on } K_{\rho_{n+1}}, \qquad \zeta_n = 0 \quad \text{on } \partial K_{\rho_n},
\]
and
\[
|\partial_{x_i}\zeta_n| \le C 2^n \rho_0^{-p^+/p_i} \omega^{(p^+ - p_i)/p_i}.
\]

Applying the energy estimate (3.1) to \((u-k_n)_-\) on \(Q_n\) gives
\[
\begin{aligned}
\sup_{t\in(-t_0,0)} &\int_{K_{\rho_n}} (u-k_n)_-^2 \zeta_n^\alpha \, dx
+ \sum_{i=1}^N \int_{Q_n} |\partial_{x_i}(u-k_n)_-|^{p_i} \zeta_n^\alpha \, dxdt \\
&\le C \sum_{i=1}^N \int_{Q_n} (u-k_n)_-^{p_i} |\partial_{x_i}\zeta_n|^{p_i} \, dxdt.
\end{aligned}
\]

On the support of the integrands, we have the sharp bound
\[
(u-k_n)_- \le k_n - \mu^- \le \frac{\omega}{2^{n_2}}.
\]
Therefore, for each \(i\),
\[
\int_{Q_n} (u-k_n)_-^{p_i} |\partial_{x_i}\zeta_n|^{p_i} \, dxdt
\le C 2^{n p_i} \rho_0^{-p^+} \omega^{p^+} 2^{-n_2 p_i} |A_n|.
\]
Since \(p_i \ge p^- := \min\{p_1,\dots,p_N\}\), we have \(2^{-n_2 p_i} \le 2^{-n_2 p^-}\). Hence
\[
\begin{aligned}
\sup_{t\in(-t_0,0)} &\int_{K_{\rho_n}} (u-k_n)_-^2 \zeta_n^\alpha \, dx
+ \sum_{i=1}^N \int_{Q_n} |\partial_{x_i}(u-k_n)_-|^{p_i} \zeta_n^\alpha \, dxdt \\
&\le C 2^{n p^+} \rho_0^{-p^+} \omega^{p^+} 2^{-n_2 p^-} |A_n|.
\end{aligned}
\]

Set \(v_n := (u-k_n)_- \zeta_n^{\alpha/2}\). Applying the anisotropic embedding Lemma 2.2 with \(l = \bar p(1+2/N)\), we obtain
\[
\int_{Q_n} |v_n|^l \, dxdt
\le C \left( 2^{n p^+} \rho_0^{-p^+} \omega^{p^+} 2^{-n_2 p^-} |A_n| \right)^{\frac{N+\bar p}{N}}.
\]

On \(A_{n+1} = \{u < k_{n+1}\}\), we have
\[
(u-k_n)_- \ge k_n - k_{n+1} \ge \frac{\omega}{2^{n_2+n+4}}.
\]
Therefore,
\[
|A_{n+1}| \left( \frac{\omega}{2^{n_2+n+4}} \right)^l
\le \int_{A_{n+1}} (u-k_n)_-^l \, dxdt
\le \int_{Q_n} |v_n|^l \, dxdt.
\]
Combining the last two estimates and solving for \(|A_{n+1}|\) gives
\[
|A_{n+1}|
\le C 2^{n (l + p^+ A)} 2^{n_2 (l - p^- A)}
\rho_0^{-p^+ A} \omega^{p^+ A - l} |A_n|^A,
\]
where
\[
A := \frac{N+\bar p}{N}.
\]

The volume of \(Q_n\) is the same as in the proof of Lemma 5.1, with \(\rho_n \sim \rho_0\). Define
\[
X_n := \frac{|A_n|}{|Q_n|}.
\]
Dividing the previous inequality by \(|Q_{n+1}| \sim |Q_n|\) and using the intrinsic cancellations ( similar to the proof of Lemma 5.1), we obtain
\[
X_{n+1} \le C 2^{n_2 \kappa} B^n X_n^{1+\alpha},
\]
where
\[
\alpha := \frac{\bar p}{N}, \qquad B := 2^{l + p^+ A},
\]
and
\[
\kappa := l - p^- A = (\bar p - p^-) + \frac{\bar p(2 - p^-)}{N}.
\]

By the assumed closeness condition \(N(\bar p - p^-) \le \bar p(p^- - 2)\), we have \(\kappa \le 0\). Hence \(2^{n_2 \kappa} \le 1\), and the recursion simplifies to
\[
X_{n+1} \le C B^n X_n^{1+\alpha}.
\]

By Lemma 2.5, \(X_n \to 0\) provided
\[
X_0 \le \nu_* := C^{-1/\alpha} B^{-1/\alpha^2}.
\]
Since \(\nu_* > 0\) is a fixed constant depending only on the data, and Lemma 5.2 gives
\[
X_0 \le \frac{C_{\log}}{n_2},
\]
we may choose \(n_2\) (which is the same as in Lemma 5.2) large enough so that
\[
\frac{C_{\log}}{n_2} \le \nu_*.
\]
Thus \(X_n \to 0\). Since
\[
k_n \downarrow \mu^- + \frac{\omega}{2^{n_2+1}},
\]
we conclude that
\[
\left| \left\{ u \le \mu^- + \frac{\omega}{2^{n_2+1}} \right\} \cap K_{\frac{\rho_0}{8}} \times (-t_0, 0) \right| = 0.
\]
Hence
\[
u > \mu^- + \frac{\omega}{2^{n_2+1}} \quad \text{a.e. in } K_{\frac{\rho_0}{8}} \times (-t_0, 0).
\]
This proves the lemma.
\end{proof}

\subsection{Second alternative}
In this subsection, we assume that \eqref{eq:second-alt} holds. Then there
exists
\[
\tau_0\in
\left[
\tau^*-\theta\rho_0^{p_+},
\tau^*-\frac{\nu_0}{2}\theta\rho_0^{p_+}
\right]
\]
such that
\begin{equation}
    \left|
    \left\{
    x\in K_{\rho_0}:
    u(x,\tau_0)>\mu^+-\frac{\omega}{2}
    \right\}
    \right|
    \le
    \left(
    \frac{1-\nu_0}{1-\frac{\nu_0}{2}}
    \right)
    |K_{\rho_0}|.
    \label{eq:tau0-second-alt}
\end{equation}
If \eqref{eq:tau0-second-alt} were false for every \(\tau_0\) in the above
interval, then for every such \(\tau_0\),
\[
\left|
\left\{
x\in K_{\rho_0}:
u(x,\tau_0)\ge \mu^-+\frac{\omega}{2}
\right\}
\right|
>
\frac{1-\nu_0}{1-\frac{\nu_0}{2}}
|K_{\rho_0}|.
\]
Integrating over
\[
\left[
\tau^*-\theta\rho_0^{p_+},
\tau^*-\frac{\nu_0}{2}\theta\rho_0^{p_+}
\right]
\]
would give
\[
\begin{aligned}
&\left|
\left\{
(x,t)\in Q^*:
u(x,t)\ge \mu^-+\frac{\omega}{2}
\right\}
\right|\\
&\ge
\left|
\left\{
(x,t)\in K_{\rho_0}\times
\left[
\tau^*-\theta\rho_0^{p_+},
\tau^*-\frac{\nu_0}{2}\theta\rho_0^{p_+}
\right]:
u(x,t)\ge \mu^-+\frac{\omega}{2}
\right\}
\right|\\
&>
\left(1-\frac{\nu_0}{2}\right)\theta\rho_0^{p_+}
\cdot
\frac{1-\nu_0}{1-\frac{\nu_0}{2}}
|K_{\rho_0}|\\
&=
(1-\nu_0)\theta\rho_0^{p_+}|K_{\rho_0}|\\
&=
(1-\nu_0)|Q^*|,
\end{aligned}
\]
contradicting \eqref{eq:second-alt}. Hence
\eqref{eq:tau0-second-alt} holds.

\begin{lem}
There exists a positive integer \(n_3>1\), depending only on the data,
such that
\begin{equation}
    \left|
    \left\{
    x\in K_{\rho_0}:
    u(x,t)>\mu^+-\frac{\omega}{2^{n_3}}
    \right\}
    \right|
    \le
    \left(
    1-\left(\frac{\nu_0}{2}\right)^2
    \right)
    |K_{\rho_0}|,
    \qquad
    \forall t\in(-\frac{a_0}{2}\rho_0^{p^+},0).
    \label{eq:second-alt-time-slices}
\end{equation}
\end{lem}

\begin{proof}
We use the intrinsic logarithmic estimate (3.16) for \(\Gamma_+\). Let
\[
k:=\mu^+-\frac{\omega}{2}
=
\frac12(\mu^++\mu^-).
\]
Since \(\underset{Q^*}{\operatorname*{ess\,sup}}~\phi\le \mu^-+\omega/2=k\), then by (5.8)
the level \(k\) is admissible for the positive logarithmic estimate.\\
Set
\[
H_k^+:=
\operatorname*{ess\,sup}_{
K_{\rho_0}\times(\tau_0,\tau^*)
}
(u-k)_+
\le
\frac{\omega}{2}.
\]
For \(c=\omega/2^{n+1}\), consider the logarithmic function
\[
\Gamma_+(u)
=
\ln\left(
\frac{H_k^+}{H_k^+-(u-k)_++c}
\right)_+ .
\]
On the set where \(u>\mu^+-\omega/2\), we have \((u-k)_+>0\), and
\[
\Gamma_+(u)\le n\ln 2,
\qquad
0\le \Gamma_+'(u)\le \frac{1}{c},
\qquad
\left|\Gamma_+'(u)\right|^{2-p_i}
\le
C\left(\frac{\omega}{2}\right)^{p_i-2}.
\]
Apply (3.16) on the nested boxes \(K_{(1-\sigma)\rho_0}\subset K_{\rho_0}\) over the time interval \((\tau_0,\tau^*)\):
\[
\begin{aligned}
&\operatorname*{ess\,sup}_{t\in(\tau_0,\tau^*)}
\int_{K_{(1-\sigma)\rho_0}}\Gamma_+^2(u)\,dx\\
&\le
\int_{K_{\rho_0}}\Gamma_+^2(u(\cdot,\tau_0))\,dx
+
C\sum_{i=1}^N
\frac{1}{(L_i(\rho_0)-L_i((1-\sigma)\rho_0))^{p_i}}
\iint_{K_{\rho_0}\times(\tau_0,\tau^*)}
\Gamma_+(\Gamma_+')^{2-p_i}\,dxdt.
\end{aligned}
\]
Here
\[
L_i(\rho)=\rho^{\frac{p_+}{p_i}}\omega^{\frac{p_i-p_+}{p_i}}.
\]
Since
\[
L_i(\rho_0)-L_i((1-\sigma)\rho_0)
\simeq
\sigma\rho_0^{\frac{p_+}{p_i}}\omega^{\frac{p_i-p_+}{p_i}},
\]
we have
\[
\frac{1}{(L_i(\rho_0)-L_i((1-\sigma)\rho_0))^{p_i}}
\le
C\sigma^{-p_i}\rho_0^{-p_+}\omega^{p_+-p_i}.
\]
The initial term is estimated using the time-slice measure from the
second alternative. By \eqref{eq:tau0-second-alt}, and since \(\Gamma_+^2\le n^2(\ln2)^2\) on \(K_{\rho_0}\), and
\(\Gamma_+=0\) whenever \(u\le \mu^+-\omega/2\), we have
\[
\begin{aligned}
\int_{K_{\rho_0}}\Gamma_+^2(u(\cdot,\tau_0))\,dx
&\le
n^2(\ln2)^2
\left|
\left\{
x\in K_{\rho_0}:\Gamma_+(u(x,\tau_0))>0
\right\}
\right| \\
&\le
n^2(\ln2)^2
\left|
\left\{
x\in K_{\rho_0}:u(x,\tau_0)>\mu^+-\frac{\omega}{2}
\right\}
\right| \\
&\le
n^2(\ln2)^2
\frac{1-\nu_0}{1-\frac{\nu_0}{2}}
|K_{\rho_0}|.
\end{aligned}
\]
For the double integral,
\[
\iint
\Gamma_+(\Gamma_+')^{2-p_i}\,dxdt
\le
C n
\left(\frac{\omega}{2}\right)^{p_i-2}
|K_{\rho_0}|(\tau^*-\tau_0).
\]
Since \(\tau^*-\tau_0\le \theta\rho_0^{p_+}\), and
\[
\theta=\left(\frac{\omega}{2}\right)^{2-p_+},
\]
we have
\[
\begin{aligned}
\frac{1}{(L_i(\rho_0)-L_i((1-\sigma)\rho_0))^{p_i}}&
\iint
\Gamma_+(\Gamma_+')^{2-p_i}\,dxdt\\
&\le
C
\sigma^{-p_i}\rho_0^{-p_+}\omega^{p_+-p_i}
\cdot
n
\left(\frac{\omega}{2}\right)^{p_i-2}
|K_{\rho_0}|
\cdot
\left(\frac{\omega}{2}\right)^{2-p_+}
\rho_0^{p_+}\\
&=
C
\sigma^{-p_i} n
|K_{\rho_0}|.
\end{aligned}
\]
Therefore
\[
\operatorname*{ess\,sup}_{t\in(\tau_0,\tau^*)}
\int_{K_{(1-\sigma)\rho_0/2}}\Gamma_+^2(u)\,dx
\le
n^2(\ln2)^2
\frac{1-\nu_0}{1-\frac{\nu_0}{2}}
|K_{\rho_0}|
+
C
\sigma^{-p_+}
n
|K_{\rho_0}|.
\]
Now fix \(t\in(\tau_0,\tau^*)\). Thereafter, define
\[
S
=
\left\{
x\in K_{(1-\sigma)\rho_0}:
u(x,t)>\mu^+-\frac{\omega}{2^{n+1}}
\right\}.
\]
Then in this set, we have
\[
\frac{H_k^+}{H_k^+-(u-k)_++c}
\ge
2^{n-1},
\]
hence
\[
\Gamma_+^2(u(x,t))
\ge
(n-1)^2(\ln2)^2.
\]
Therefore,
\[
(n-1)^2(\ln2)^2 |S|
\le
n^2(\ln2)^2
\frac{1-\nu_0}{1-\frac{\nu_0}{2}}
|K_{\rho_0}|
+
C
\sigma^{-p_+}
n
|K_{\rho_0}|.
\]
On the other hand, for all $t\in(\tau_{0},\tau^{*})$, we get
\begin{equation}
    \begin{split}
     \biggl|\biggl\{x\in K_{\rho_{0}}&:~u>\mu^{+}-\frac{\omega}{2^{n+1}}\biggr\}\biggr| \leq \left|S\right|+N\sigma|K_{\rho_{0}}|\\
     &\leq \left\{ \left(\frac{n}{n-1}\right)^{2}\left(\frac{1-\nu_{0}}{1-\frac{\nu_{0}}{2}}\right) +\frac{C}{n\sigma^{p^{+}}} +N\sigma  \right\}|K_{\rho_{0}}|\\
     &\leq \left( 1-\left( \frac{\nu_{0}}{2}\right)^{2}  \right)|K_{\rho_{0}}|,
    \end{split}
\end{equation}
where we took $\sigma$ small enough such that $N\sigma\leq \frac{3}{8}\nu_{0}^{2}$ and $n$ so large that $ \left(\frac{n}{n-1}\right)^{2}\leq (1-\frac{\nu_{0}}{2})(1+\nu_{0})$ and $\frac{C}{n\sigma^{p^{+}}}\leq \frac{3}{8}\nu_{0}^{2}$. Finally, recalling that \(\tau_0\in[\tau^*-\theta\rho_0^{p^+},\tau^*-\frac{\nu_0}{2}\theta\rho_0^{p^+}]\) and choosing \(\lambda\) such that \(2^{(\lambda-1)(p^+-2)}\ge2\), we get the lemma for all \(t\in(-\frac{a_0}{2}\rho_0^{p^+},0)\).
\end{proof}
Now, we are going to use the result of Lemma 5.4 to get that within the cylinder \(\mathcal{Q}_{\frac{\theta}{2},\omega}(\rho_0)\), the set of points where \(u\) is close to its supremum has an arbitrarily small measure.

\begin{lem}
    For \(\tilde{\nu}_{1}\in(0,1)\), there exists an integer \(\lambda\geq n_3\) depending on the data such that
    \begin{equation}
        \left|
        \left\{
            (x,t)\in \mathcal{Q}_{\frac{a_0}{2},\omega}(\rho_0):
            u(x,t)>\mu^+-\frac{\omega}{2^{\lambda}}
        \right\}
        \right|
        \leq
        \tilde{\nu}_{1}
        \left|
        \mathcal{Q}_{\frac{a_0}{2},\omega}(\rho_0)
        \right|.
        \label{eq:close-to-supremum-small}
    \end{equation}
\end{lem}
\begin{proof}
We begin by taking $k=\mu^{+}-\frac{\omega}{2^{n}}\geq\mu^{-}+\frac{\omega}{2}\geq\underset{\mathcal{Q}_{a_0,\omega}(\rho_0)}{ess\sup}~\phi$ for $n_{3}\leq n\leq\lambda$. Therefore, we can apply (3.2) for $(u-k)_{+}$ where we take $0\leq\xi(x,t)\leq1$ as a smooth cutoff function satisfying  
\begin{equation*}
    \begin{cases}
        \xi=1~~\text{in}~\mathcal{Q}_{\frac{a_0}{2},\omega}(\rho_0),~\xi=0~\text{on }~\partial_{p}\mathcal{Q}_{a_0,\omega}(2\rho_0),
        &\\
        \left|\frac{\partial \xi}{\partial x_{i}} \right|\leq C\rho_{0}^{-p_{+}/p_i}\omega^{(p_+ - p_i)/p_i}~\text{for }~i=1,..,N,~0<\frac{\partial \xi}{\partial t} \leq Ca_0^{-1}\rho_{0}^{-p_+},
    \end{cases}
\end{equation*}
such that for $n\leq\lambda$
\begin{equation}
    \sum_{i=1}^{N}\int_{\mathcal{Q}_{\frac{a_0}{2},\omega}(\rho_0)}\left| \frac{\partial}{\partial x_{i}}(u-k)_{+}   \right|^{p_{i}}~dxdt\leq \frac{C}{\rho_{0}^{p_+}}\left(\frac{\omega}{2^{n}}   \right)^{p_+}|\mathcal{Q}_{\frac{a_0}{2},\omega}(\rho_0)|,
\end{equation}
where we used the same method and similar assumptions as the ones we used to get (5.13). Now, for $n\leq \lambda$ we define the following sets
\begin{equation*}
    G_{n}(t)=\{x\in K_{\rho_{0}}: u(x,t)>\mu^{+}-\frac{\omega}{2^{n}}\}, \qquad 
    G_{n}=\int_{-\frac{a_0\rho_{0}^{p_+}}{2}}^{0}G_{n}(t)\,dt,
\end{equation*}
and 
\begin{equation*}
    K_{\rho_{0}}\setminus G_{n}(t)=\{x\in K_{\rho_{0}}: u(x,t)\leq\mu^{+}-\frac{\omega}{2^{n}}\}.
\end{equation*}
Also, for all $t\in(-\frac{a_0\rho_{0}^{p_+}}{2},0)$ we define the truncation function
\begin{equation}
   \gamma_{n}(x,t)= \begin{cases}
        0 & \text{for } u<\mu^{+}-\dfrac{\omega}{2^{n}},\\[2mm]
        u-\left(\mu^{+}-\dfrac{\omega}{2^{n}}\right) & \text{for } \mu^{+}-\dfrac{\omega}{2^{n}}\leq u<\mu^{+}-\dfrac{\omega}{2^{n+1}},\\[2mm]
        \dfrac{\omega}{2^{n+1}} & \text{for } u\geq \mu^{+}-\dfrac{\omega}{2^{n+1}}.
  \end{cases}
\end{equation}
By construction, $\gamma_{n}$ vanishes on the set $K_{\rho_{0}}\setminus G_{n}(t)$.

Now, since $K_{\rho_{0}}$ is a Cartesian product of intervals, for any $x=(x_1,\dots,x_N)\in G_{n}(t)$ and any $y=(y_1,\dots,y_N)\in K_{\rho_{0}}\setminus G_{n}(t)$, the coordinate-wise polygonal path connecting $y$ to $x$ remains entirely inside $K_{\rho_{0}}$. For example, we take
\[
\pi_N = x,\quad \pi_{N-1} = (x_1,\dots,x_{N-1},y_N),\quad \dots,\quad \pi_0 = y.
\]
Each segment of this path lies in $K_{\rho_0}$ because each coordinate stays within its respective interval. Therefore, by the fundamental theorem of calculus,
\begin{equation}
     \begin{aligned}
\gamma_{n}(x,t)
&=[\gamma_{n}(\pi_N,t)-\gamma_{n}(\pi_{N-1},t)]+\cdots+[\gamma_{n}(\pi_1,t)-\gamma_{n}(\pi_0,t)] \\
&=\int_{y_N}^{x_N}\partial_{x_N}\gamma_{n}(x_1,\dots,x_{N-1},\zeta,t)\,d\zeta
+\cdots
+\int_{y_1}^{x_1}\partial_{x_1}\gamma_{n}(\zeta,x_2,\dots,x_N,t)\,d\zeta \\
&\leq \sum_{i=1}^{N}\int_{-L_i}^{L_i}
\left|\partial_{x_i}\gamma_{n}(x_1,\dots,\underset{i\text{-th}}{\zeta},\dots,x_N,t)\right|\,d\zeta,
\end{aligned}
\end{equation}
where
\[
L_i := \rho_0^{p_+/p_i}\omega^{(p_i-p_+)/p_i}
\]
is the half-length of $K_{\rho_0}$ in the $x_i$-direction.

Now double integrate this inequality in $dx$ over $x\in G_n(t)$ and in $dy$ over $y\in K_{\rho_0}\setminus G_n(t)$. The left-hand side becomes
\[
\int_{G_n(t)}\int_{K_{\rho_0}\setminus G_n(t)} \gamma_n(x,t)\,dy\,dx
= |K_{\rho_0}\setminus G_n(t)| \int_{G_n(t)} \gamma_n(x,t)\,dx.
\]
By Lemma 5.4, for every $t\in(-\frac{a_0}{2}\rho_0^{p_+},0)$,
\[
|K_{\rho_0}\setminus G_n(t)| \ge \left(\frac{\nu_0}{2}\right)^2 |K_{\rho_0}|.
\]
Thus the left-hand side is bounded below by
\[
\left(\frac{\nu_0}{2}\right)^2 |K_{\rho_0}| \int_{G_n(t)} \gamma_n(x,t)\,dx.
\]

For the right-hand side, after double integration, the factor $2L_i$ appears (the length of the full interval in the $x_i$-direction), where
\[
L_i := \rho_0^{p_+/p_i}\omega^{(p_i-p_+)/p_i}
\]
is the half-length of $K_{\rho_0}$ in the $x_i$-direction. Enlarging the domains to the full $K_{\rho_0}$, we obtain
\[
\text{RHS} \le |K_{\rho_0}| \sum_{i=1}^{N} 2L_i \int_{K_{\rho_0}} |\partial_{x_i}\gamma_n(x,t)|\,dx.
\]
Combining this with the lower bound from Lemma 5.4 gives
\[
\left(\frac{\nu_0}{2}\right)^2 |K_{\rho_0}| \int_{G_n(t)} \gamma_n(x,t)\,dx
\le
|K_{\rho_0}| \sum_{i=1}^{N} 2L_i \int_{K_{\rho_0}} |\partial_{x_i}\gamma_n(x,t)|\,dx.
\]
Since $\gamma_n=0$ on $K_{\rho_0}\setminus G_n(t)$, we have
\[
\int_{G_n(t)} \gamma_n(x,t)\,dx = \int_{K_{\rho_0}} \gamma_n(x,t)\,dx.
\]
Therefore,
\begin{equation}
    \left(\frac{\nu_0}{2} \right)^{2}|K_{\rho_{0}}|\int_{K_{\rho_{0}}}\gamma_n(x,t)\,dx
    \leq |K_{\rho_{0}}| \sum_{i=1}^{N} 2L_i \int_{K_{\rho_{0}}}\left|\partial_{x_i}\gamma_n(x,t)\right|\,dx.
\end{equation}
Since $\gamma_n = \omega/2^{n+1}$ on $G_{n+1}(t)$, we have
\[
\int_{K_{\rho_0}} \gamma_n(x,t)\,dx \ge \frac{\omega}{2^{n+1}} |G_{n+1}(t)|.
\]
Also, $\partial_{x_i}\gamma_n = \partial_{x_i}u$ on $G_n(t)\setminus G_{n+1}(t)$ and zero elsewhere. Substituting these into (5.26) yields
\begin{equation}
   \frac{\omega}{2^{n+1}}|G_{n+1}(t)|
   \leq \frac{C}{\nu_{0}^{2}}\sum_{i=1}^{N} L_i
   \int_{G_{n}(t)\setminus G_{n+1}(t)} |\partial_{x_i}u|\,dx,
\end{equation}
Then, by integrating (5.27) over $t\in(-\frac{a_0}{2}\rho_0^{p_+},0)$ and using (5.23), we obtain
\begin{equation}
     \begin{split}
|G_{n+1}|&\leq \frac{C}{\nu_0^2}\frac{2^{n+1}}{\omega}\sum_{i=1}^{N} L_i
\left(\int_{D_n}\left|\frac{\partial u}{\partial x_{i}}   \right|^{p^-}~dxdt   \right)^{\frac{1}{p^-}}|D_n|^{\frac{p^- -1}{p^-}}\\
&\leq \frac{C}{\nu_0^2}
\sum_{i=1}^{N} 2^{n+1 - n p_+/p_i}
\left|\mathcal{Q}_{\frac{a_0}{2},\omega}(\rho_0)\right|^{\frac{1}{p^-}}
|D_n|^{\frac{p^- -1}{p^-}}\\
&\leq \frac{C}{\nu_0^2}
\left|\mathcal{Q}_{\frac{a_0}{2},\omega}(\rho_0)\right|^{\frac{1}{p^-}}
|D_n|^{\frac{p^- -1}{p^-}},
     \end{split}
\end{equation}
where \(D_n := G_n \setminus G_{n+1}\), in the last inequality we used \(2^{n+1 - n p_+/p_i} \le 2\) (since \(p_i \le p_+\)) and absorbed the factor 2 into the constant \(C\).

By raising (5.28) to the power $\frac{p^{-}}{p^{-}-1}$ and summing it up from $n=n_{3},n_{3}+1,\dots,\lambda-1$, we obtain
\begin{equation}
    \begin{split}
\sum_{n=n_{3}}^{\lambda-1}|G_{n+1}|^{\frac{p^{-}}{p^{-}-1}}
&\leq C(\nu_{0})^{-\frac{2p^{-}}{p^{-}-1}}
\left|\mathcal{Q}_{\frac{a_0}{2},\omega}(\rho_0)\right|^{\frac{1}{p^{-}-1}}
\sum_{n=n_{3}}^{\lambda-1}|G_{n}-G_{n+1}|.
    \end{split}
\end{equation}

Next, since the sets $G_n$ are decreasing (i.e., $G_{\lambda} \subset G_{n+1}$ for all $n \le \lambda-1$), we have $|G_{\lambda}| \le |G_{n+1}|$. Therefore,
\[
\sum_{n=n_{3}}^{\lambda-1}|G_{n+1}|^{\frac{p^{-}}{p^{-}-1}}
\geq (\lambda-n_3)|G_{\lambda}|^{\frac{p^{-}}{p^{-}-1}}.
\]
Also, the sets $G_n \setminus G_{n+1}$ are pairwise disjoint and their union is contained in $\mathcal{Q}_{\frac{a_0}{2},\omega}(\rho_0)$, so
\[
\sum_{n=n_{3}}^{\lambda-1}|G_{n}-G_{n+1}|
= |G_{n_3} \setminus G_{\lambda}|
\le \left|\mathcal{Q}_{\frac{a_0}{2},\omega}(\rho_0)\right|.
\]
Substituting these two estimates into (5.29) gives
\[
(\lambda-n_3)|G_{\lambda}|^{\frac{p^{-}}{p^{-}-1}}
\le
C(\nu_{0})^{-\frac{2p^{-}}{p^{-}-1}}
\left|\mathcal{Q}_{\frac{a_0}{2},\omega}(\rho_0)\right|^{\frac{1}{p^{-}-1}}
\left|\mathcal{Q}_{\frac{a_0}{2},\omega}(\rho_0)\right|.
\]
Since
\[
\frac{1}{p^{-}-1} + 1 = \frac{p^{-}}{p^{-}-1},
\]
we get
\[
(\lambda-n_3)|G_{\lambda}|^{\frac{p^{-}}{p^{-}-1}}
\le
C(\nu_{0})^{-\frac{2p^{-}}{p^{-}-1}}
\left|\mathcal{Q}_{\frac{a_0}{2},\omega}(\rho_0)\right|^{\frac{p^{-}}{p^{-}-1}}.
\]
Taking the $(p^{-}-1)/p^{-}$-th power yields
\begin{equation}
    |G_{\lambda}|
    \leq \frac{C}{(\lambda-n_3)^{\frac{p^{-}-1}{p^{-}}}}
    (\nu_{0})^{-2}
    \left|\mathcal{Q}_{\frac{a_0}{2},\omega}(\rho_0)\right|.
\end{equation}
Hence, by choosing $\lambda$ large enough such that
\[
\frac{C}{(\lambda-n_3)^{\frac{p^{-}-1}{p^{-}}}}(\nu_{0})^{-2}
\leq \widetilde{\nu}_{1} < 1,
\]
we get the desired result
\[
|G_{\lambda}|
\le \widetilde{\nu}_{1}
\left|\mathcal{Q}_{\frac{a_0}{2},\omega}(\rho_0)\right|.
\]
\end{proof}

\begin{lem}
For $\widetilde{\nu}_1\in(0,1)$, the choice of $\lambda$ from Lemma 5.5 can be made so that
\begin{equation*}
    u \le \mu^{+} - \frac{\omega}{2^{\lambda+1}} \quad \text{a.e. in } 
    \mathcal{Q}_{\frac{a_0}{2},\omega}\!\left(\frac{\rho_0}{2}\right).
\end{equation*}
\end{lem}
\begin{proof}
We begin our proof by taking decreasing sequences
\begin{equation*}
    \rho_{n}=\frac{\rho_{0}}{2}+\frac{\rho_{0}}{2^{n+1}},\qquad
    k_{n}=\mu^{+}-\frac{\omega}{2^{\lambda+1}}-\frac{\omega}{2^{\lambda+1+n}},\qquad n=0,1,\dots
\end{equation*}
Therefore, since $k_{n}\geq\mu^{+}-\frac{\omega}{2^{\lambda}}\geq \mu^{-}+\frac{\omega}{2}\geq\underset{\mathcal{Q}_{a_0,\omega}(\rho_0)}{\operatorname{ess}\,\sup}~\phi$, which is guaranteed by (5.8), we can apply (3.2) to $(u-k_{n})_{+}$ on the intrinsic cylinders $\mathcal{Q}_{\frac{a_0}{2},\omega}(\rho):=K_\rho \times \left(-\frac{a_0}{2}\rho^{p_+},0\right)$. Next, we choose a smooth cutoff $0\le \xi_n(x,t)\le 1$ such that
\begin{equation*}
    \begin{cases}
        \xi_n \equiv 1 ~~\text{in } \mathcal{Q}_{\frac{a_0}{2},\omega}(\rho_{n+1}), \qquad
        \xi_n \equiv 0 ~~\text{on the parabolic boundary of } \mathcal{Q}_{\frac{a_0}{2},\omega}(\rho_n),
        \\[2mm]
        \left|\partial_{x_i}\xi_n\right| \le \dfrac{C\,2^n}{L_i(\rho_n) - L_i(\rho_{n+1})}
        \quad (i=1,\dots,N), \\[2mm]
        \left|\partial_t \xi_n\right| \le \dfrac{C\,2^n}{a_0\left(\rho_n^{p_+} - \rho_{n+1}^{p_+}\right)},
    \end{cases}
\end{equation*}
where $L_i(\rho):=\rho^{p_+/p_i}\omega^{(p_i-p_+)/p_i}$ is the half‑length of $K_\rho$ in the $x_i$-direction. The rest of the proof is similar to the proof of Lemma 5.3 by using Lemma 5.5.
\end{proof}
\subsection{The Recursive Argument}
Based on our previous findings, we can conclude that the oscillation of $u$ is reduced in both alternatives.
\begin{cor}
 There exists $\sigma\in(0,1)$ depending on the data such that
 \begin{equation*}
\underset{\mathcal{Q}_{\theta,\omega}(\frac{\rho_0}{8})}{ess~osc}~u\leq\sigma\omega.
 \end{equation*}
\end{cor}
\begin{proof}
    From Lemma 5.3 we get that
    \begin{equation}
    \begin{split}
        \underset{\mathcal{Q}_{\theta,\omega}(\frac{\rho_0}{8})}{ess~osc}~u&\leq \underset{K_{\frac{\rho_0}{2}}\times(-t_0,0)}{ess\sup}~u-\underset{K_{\frac{\rho_0}{2}}\times(-t_0,0)}{ess\inf}~u\\
        &\leq \mu^{+}-\mu^{-}-\frac{\omega}{2^{n_{2}+1}}=\left(1-\frac{1}{2^{n_{2}+1}}  \right)\omega=\sigma_{1}\omega.
        \end{split}
    \end{equation}
    Next, from Lemma 5.6 we get also that
     \begin{equation}
        \underset{\mathcal{Q}_{\theta,\omega}(\frac{\rho_0}{8})}{ess~osc}~u\leq \sigma_{2}\omega,
    \end{equation}
    where $\sigma_{2}=\left( 1-\frac{1}{2^{\lambda+1}} \right)$. Hence, from (5.31) and (5.32) we get the desired result for $\sigma=\max\{\sigma_{1},\sigma_{2}\}.$
\end{proof}
Consequently, we obtain the following recursive result
\begin{prop}
    For $\sigma\in(0,1)$ and by letting
    \begin{equation*}
        \omega_{1}=\max\{\sigma\omega, 2~\underset{\mathcal{Q}_{a_0,\omega}(\rho_0)}{ess~osc}~\phi \},
    \end{equation*}
    there exists a positive constant $\gamma$ depending on the data such that
    \begin{equation*}
        \gamma=\sigma^{\frac{p^{+}-2}{p^{+}}}2^{\frac{(\lambda-1)(2-p^{+})}{p^{+}}-3},
    \end{equation*}
    and by constructing the cylinder
    \begin{equation*}
    Q_{1}=\mathcal{Q}_{a_1,\omega_1}(\rho_1),~a_{1}=\left( \frac{\omega_{1}}{2^{\lambda}} \right)^{2-p^{+}},~\rho_{1}=\gamma\rho_{0},    
    \end{equation*}
 we have that
\begin{equation*}
    \underset{Q_{1}}{ess~osc}~u\leq\omega_{1},~~\text{and}~Q_{1}\subset \mathcal{Q}_{a_0,\omega}(\rho_0).
\end{equation*}
\end{prop}
\begin{proof}
By construction, from (5.7) we have that
\begin{equation*}
      \underset{\mathcal{Q}_{a_0,\omega}(\rho_0)}{ess~osc}~u\leq\omega.
\end{equation*}
Next, since $\sigma\omega\leq\omega_{1}$, we find that 
\begin{equation*}
    \begin{split}
\theta\left( \frac{\rho_{0}}{8} \right)^{p^{+}}&=\left( \frac{\omega}{2} \right)^{2-p^{+}}\left( \frac{2^{\lambda}}{\omega_{1}} \right)^{2-p^{+}}
\left( \frac{\omega_{1}}{2^{\lambda}} \right)^{2-p^{+}}\frac{\rho_{0}^{p^{+}}}{2^{3p^{+}}}\\
&=\left( \frac{\omega}{\omega_{1}} \right)^{2-p^{+}}2^{(\lambda-1)(2-p^{+})-3p^{+}}\theta_{1}\rho_{0}^{p^{+}}\\
&\geq \sigma^{p^{+}-2}2^{(\lambda-1)(2-p^{+})-3p^{+}}\theta_{1}\rho_{0}^{p^{+}}\\
&=\rho_{1}^{p^{+}}\theta_{1}.
    \end{split}
\end{equation*}
For spatial inclusion, by simple computation we get $\gamma \le \frac{1}{8} \sigma^{(p_+ - p_i)/p_+}$ for each $i$, which implies that  $K_{\rho_1} \subset K_{\rho_0/8}$. Therefore, by using Corollary 5.7 and the definition of $\omega_{1}$, we get that
\begin{equation*}
     \underset{Q_{1}}{ess~osc}~u\leq\sigma\omega\leq\omega_{1}.
\end{equation*}
\end{proof}
\subsection{Proof of Theorem 1.1.}
We begin by defining the following sequences of parameters and sequences such that for $n=1,2,..$
\begin{equation}
   \begin{cases}
\rho_{n}=\gamma\rho_{n-1},~\omega_{n}=\max\{\sigma\omega_{n-1},~2~\underset{Q_{n-1}}{ess~osc}~\phi\},~a_{n}=\left(\frac{\omega_{n}}{2^{\lambda}} \right)^{2-p^{+}},\\
 \gamma=\sigma^{\frac{p^{+}-2}{p^{+}}}2^{\frac{(\lambda-1)(2-p^{+})}{p^{+}}-3}\in(0,1),~Q_{n}=\mathcal{Q}_{a_n,\omega_n}(\rho_n),\\
 \mu^{-}_{n}=\underset{Q_{n}}{ess\inf}~u,~~\text{and}~\mu^{+}_{n}=\mu_{n}^{-}+\omega_{n}.
   \end{cases}
\end{equation}
By proposition 5.8, we have that
\begin{equation}
    Q_{n}\subset Q_{n-1}~~~~\text{and}~~~~\underset{Q_{n}}{ess~osc}~u\leq\omega_{n}.
\end{equation}
Indeed, for $n=0$ the result is assured by (5.7) for $\omega_{0}:=\omega$. Next, we assume that (5.34) is true for $n$ and we will prove it for (n+1). Therefore, we have $\mu^{-}_{n}=\underset{Q_{n}}{ess\inf}~u$ and
\begin{align}
   &\mu^{+}_{n}=\mu_{n}^{-}+\omega_{n}\geq \mu_{n}^{-}+\underset{Q_{n}}{ess~osc}~u=\underset{Q_{n}}{ess\sup}~u,\\
   &\underset{Q_{n}}{ess\sup}~\phi=\underset{Q_{n}}{ess\inf}~\phi+\underset{Q_{n}}{ess~osc}~\phi \leq \underset{Q_{n}}{ess\inf}~u+\underset{Q_{n-1}}{ess~osc}~\phi\leq\mu_{n}^{-}+\frac{1}{2}\omega_{n},\\
   & \text{and},~\underset{Q_{n}}{ess~osc}~\phi\leq \underset{Q_{n-1}}{ess~osc}~\phi\leq \frac{1}{2}\omega_{n}.
\end{align}
Consequently, we can apply Proposition 5.8 to obtain that
\begin{equation*}
    Q_{n+1}\subset Q_{n},~~\text{and}~~\underset{Q_{n+1}}{ess~osc}~u\leq \omega_{n+1}.
\end{equation*}
Next, we define 
\begin{equation}
    r_{n}=\min\{a^{\frac{1}{p^{+}}}_0,~\underset{1\leq i\leq N}{\min} \omega_{0}^{\frac{p_i-p^+}{p_i}}\}\rho^{\frac{p^+}{p^-}}_{n}.
\end{equation}
By induction from (5.33), we obtain
\[
\omega_n \le \omega,
\]
and consequently, we have
\begin{equation}
    Q_{r_{n}}=(-r_{n},r_{n})^N\times(-r_n^{p^{+}},0)\subset Q_{n}~~\text{for}~n\in \mathbb{N},
\end{equation}
where $a_0$ is defined in \eqref{eq:a0}. Therefore, for all $n\in\mathbb{N}$ we have
\begin{equation}
    \begin{split}
\underset{Q_{r_{n}}}{ess~osc}~u\leq \underset{Q_{n}}{ess~osc}~u&\leq \omega_{n}\leq \sigma\omega+2~\underset{Q_{n-1}}{ess~osc}~\phi\\
&\leq \sigma^{n}\omega+2\sum_{j=0}^{n-1}\sigma^{j}~\underset{Q_{n-1-j}}{ess~osc}~\phi.
\end{split}
\end{equation}
Now, we are going to simplify the last term on the right-hand side of (5.40) such that
\begin{equation}
\begin{split}
\underset{Q_{n-1-j}}{\operatorname{ess\,osc}}~\phi 
&\leq C[\phi]_{0;\beta,\frac{\beta}{2}}\left( \sum_{i=1}^{N}\left(\rho_{n-1-j}^{\frac{p^{+}}{p_{i}}}\omega_{n-1-j}^{\frac{p_{i}-p^{+}}{p_{i}}}\right)^{\beta} + \left(a_{n-1-j}\rho_{n-1-j}^{p^{+}}\right)^{\frac{\beta}{2}}\right)\\
&\leq C[\phi]_{0;\beta,\frac{\beta}{2}}\left( \sum_{i=1}^{N}\rho_{n-1-j}^{\beta\frac{p^{+}}{p_{i}}}\omega_{n-1-j}^{\beta\frac{p_{i}-p^{+}}{p_{i}}} + \left(a_{n-1-j}\rho_{n-1-j}^{p^{+}}\right)^{\frac{\beta}{p^{+}}}\right)\\
&\leq C[\phi]_{0;\beta,\frac{\beta}{2}}\left( \sum_{i=1}^{N}\omega_{n-1-j}^{\beta\frac{p_{i}-p^{+}}{p_{i}}} + \left( \frac{\omega_{n-1-j}}{2^{\lambda}} \right)^{\beta\frac{2-p^{+}}{p^{+}}} \right)\rho_{n-1-j}^{\beta}\\
&\leq C[\phi]_{0;\beta,\frac{\beta}{2}}\left( \sum_{i=1}^{N}\left( \sigma^{n-1-j}\omega \right)^{\beta\frac{p_{i}-p^{+}}{p_{i}}} + \left( \sigma^{n-1-j}\omega \right)^{\beta\frac{2-p^{+}}{p^{+}}} \right)\rho_{n-1-j}^{\beta}\\
&\leq C[\phi]_{0;\beta,\frac{\beta}{2}}\left(1 + \omega^{\beta\frac{p^{-}-p^{+}}{p^{-}}} + \omega^{\beta\frac{2-p^{+}}{p^{+}}} \right)\sigma^{\frac{\beta(2-p^{+})(n-1-j)}{p^{+}}}\rho_{n-1-j}^{\beta}\\
&\leq C[\phi]_{0;\beta,\frac{\beta}{2}}\left(1 + \omega^{\beta\frac{p^{-}-p^{+}}{p^{-}}} + \omega^{\beta\frac{2-p^{+}}{p^{+}}} \right)\left(\frac{\delta_{0}}{8}\right)^{\beta(n-1-j)}\rho_{0}^{\beta}
\end{split}
\end{equation}
where $\delta_0 = 2^{\frac{(\lambda-1)(2-p^{+})}{p^{+}}}$. We used \eqref{eq:rho0} for the second inequality, the condition $\rho_0<1$ for the third, and $\sigma^n \omega \leq \omega_n$ for the remaining estimate. Therefore,
\begin{equation}
\begin{split}
\underset{Q_{r_{n}}}{\operatorname{ess\,osc}}~u 
&\leq \sigma^{n}\omega + C[\phi]_{0;\beta,\frac{\beta}{2}}\left(1 + \omega^{\beta\frac{p^{-}-p^{+}}{p^{-}}} + \omega^{\beta\frac{2-p^{+}}{p^{+}}}\right)\sum_{j=0}^{n-1}\sigma^{j}\left(\frac{\delta_{0}}{8}\right)^{\beta(n-1-j)}\rho_{0}^{\beta}\\
&\leq \sigma^{n}\omega + C[\phi]_{0;\beta,\frac{\beta}{2}}\left(1 + \omega^{\beta\frac{p^{-}-p^{+}}{p^{-}}} + \omega^{\beta\frac{2-p^{+}}{p^{+}}}\right)n\delta_{1}^{n-1}\rho_{0}^{\beta}\\
&\leq \sigma^{n}\omega + C[\phi]_{0;\beta,\frac{\beta}{2}}\left(1 + \omega^{\beta\frac{p^{-}-p^{+}}{p^{-}}} + \omega^{\beta\frac{2-p^{+}}{p^{+}}}\right)\sqrt{\delta_{1}^{n}}\rho_{0}^{\beta}
\end{split}
\end{equation}
where $\delta_{1}=\max\{\sigma,~\left(\frac{\delta_{0}}{8} \right)^{\beta}  \}$ and we use the fact that $n\delta_{1}^{n-1}\leq C(\delta_{1})\sqrt{\delta_{1}^{n}}$. Next, defining the modified H\"older exponent $\alpha = \vartheta \frac{p^{-}}{p^{+}}$ at the outset, we set
\begin{equation}
    \vartheta=\min\left\{ \frac{\ln(\delta_{1})}{2\ln(\gamma)}, \, \frac{p^{+}\beta}{p^{+}+\beta(p^{+}-2)}, \, \frac{p^{-}\beta}{p^{-}+\beta(p^{+}-p^{-})} \right\},
\end{equation}
and note that $\vartheta\leq \frac{\ln(\sigma)}{\ln(\gamma)}$. Then, (5.42) becomes
\begin{equation}
    \underset{Q_{r_{n}}}{\operatorname{ess\,osc}}~u \leq \gamma^{n\vartheta}\omega + C[\phi]_{0;\beta,\frac{\beta}{2}}\gamma^{n\vartheta}
    \left(1 + \omega^{\beta\frac{p^{-}-p^{+}}{p^{-}}} + \omega^{\beta\frac{2-p^{+}}{p^{+}}} \right)\rho_{0}^{\beta}.
\end{equation}
From (5.3), we have
\begin{equation}
    2^{\lambda}\rho_{0}^{\vartheta} \leq \omega \leq C[\phi]_{0;\beta,\frac{\beta}{2}}R^{-\beta}\rho_{0}^{\vartheta}.
\end{equation}
Also, from the definition $r_{n} = C_{0}\rho_{n}^{\frac{p^{+}}{p^{-}}}$, where $C_{0} = \min\left\{a_{0}^{\frac{1}{p^{+}}}, \, \min_{1\le i\le N} \omega_{0}^{\frac{p_{i}-p^{+}}{p_{i}}}\right\}$, the local boundedness of $u$ ($\omega_{0} \le 2\|u\|_{L^{\infty}_{\mathrm{loc}}}$) guarantees $C_{0}^{-1} \le C(p^{\pm}, a_{0}, \|u\|_{L^{\infty}_{\mathrm{loc}}})$, giving
\begin{equation}
    \gamma^{n\vartheta} = \left(\frac{\rho_{n}}{\rho_{0}}\right)^{\vartheta} = \left( \frac{r_{n}}{C_{0}\rho_{0}^{\frac{p^{+}}{p^{-}}}} \right)^{\vartheta\frac{p^{-}}{p^{+}}} = C_{0}^{-\alpha} \frac{r_{n}^{\alpha}}{\rho_{0}^{\vartheta}} \le C(p^{\pm}, \|u\|_{L^{\infty}_{\mathrm{loc}}}) \, \frac{r_{n}^{\alpha}}{\rho_{0}^{\vartheta}}.
\end{equation}
Then, substituting (5.45) and (5.46) into (4.41), we obtain
\begin{equation}
\begin{split}
\underset{Q_{r_{n}}}{\operatorname{ess\,osc}}~u &\leq C(p^{\pm},\|u\|_{L^{\infty}_{\mathrm{loc}}},[\phi]_{0;\beta,\frac{\beta}{2}})\Bigg\{ \frac{r_{n}^{\alpha}}{R^{\beta}}
+ r_{n}^{\alpha}\rho_{0}^{\beta-\vartheta}
+ r_{n}^{\alpha}\rho_{0}^{\beta - \vartheta + \vartheta\beta\frac{p^{-}-p^{+}}{p^{-}}}
+ r_{n}^{\alpha}\rho_{0}^{\beta - \vartheta + \vartheta\beta\frac{2-p^{+}}{p^{+}}} \Bigg\}\\
&\leq C(p^{\pm},\|u\|_{L^{\infty}_{\mathrm{loc}}},[\phi]_{0;\beta,\frac{\beta}{2}})\frac{r_{n}^{\alpha}}{R^{\beta}},
\end{split}
\end{equation}
where we used the fact that $\rho_0 \leq R < 1$, $\alpha = \vartheta \frac{p^{-}}{p^{+}}$, and the choice of $\vartheta$ ensures all net exponents of $\rho_0$ are non-negative. Indeed, the three exponents are
\[
\beta-\vartheta \ge 0, \qquad 
\beta - \vartheta + \vartheta\beta\frac{p^{-}-p^{+}}{p^{-}} \ge 0, \qquad 
\beta - \vartheta + \vartheta\beta\frac{2-p^{+}}{p^{+}} \ge 0,
\]
which are precisely equivalent to the three bounds in the definition of $\vartheta$.

Thereafter, we extend estimate (5.47) to all continuous radii $r \in (0, R]$. We consider three distinct cases:

Case 1: Let $r \in (0, r_{0})$, then we choose $n \in \mathbb{N}_{0}$ such that $r_{n+1} \leq r \leq r_{n}$. Since $Q_{r} \subset Q_{r_{n}}$, the monotonicity of the essential oscillation and the scaling relation $r_{n+1}^{\alpha} = \gamma^{\vartheta} r_{n}^{\alpha}$ yield
\begin{equation}
    \underset{Q_{r}}{\operatorname{ess\,osc}}~u \leq \underset{Q_{r_{n}}}{\operatorname{ess\,osc}}~u \leq C\frac{r_{n}^{\alpha}}{R^{\beta}} = \frac{C}{\gamma^{\vartheta}} \frac{r_{n+1}^{\alpha}}{R^{\beta}} \leq \frac{C}{\gamma^{\vartheta}} \frac{r^{\alpha}}{R^{\beta}} \leq C \frac{r^{\alpha}}{R^{\beta}}.
\end{equation}

Case 2: Let $r \in [r_{0}, \rho_{0}]$, then since $r \le \rho_{0} < 1$ and $p^{+} > 2$, we have $r^{p^{+}} \le \rho_{0}^{2}$, which implies the inclusion $Q_{r} \subset Q(\rho_{0}^{2}, \rho_{0})$. From \eqref{eq:rho0}, the definition of $H(\rho_{0})$, and the scaling relation $r_{0}^{\alpha} = C_{0}^{\alpha}\rho_{0}^{\vartheta}$, we get
\begin{equation}
    \underset{Q_{r}}{\operatorname{ess\,osc}}~u \leq \underset{Q(\rho_{0}^{2},\rho_{0})}{\operatorname{ess\,osc}}~u \leq C R^{-\beta} H(\rho_{0}) \leq C R^{-\beta} \rho_{0}^{\vartheta} = C R^{-\beta} C_{0}^{-\alpha} r_{0}^{\alpha} \leq C \frac{r^{\alpha}}{R^{\beta}}.
\end{equation}

Case 3: Let $r \in [\rho_{0}, R]$, then since $r \le R < 1$ and $p^+ > 2$, we have $r^{p^+} \le r^2$, which gives the inclusion $Q_{r} \subset Q(r^{2}, r)$. Applying estimate \eqref{eq:rho0} with $\rho = r$, the definition of $H(r)$, we estimate
\begin{equation}
    \underset{Q_{r}}{\operatorname{ess\,osc}}~u \leq C \frac{r^{\alpha}}{R^{\beta}}.
\end{equation}
Hence, combining all three cases proves local H\"older continuity over $Q_{r}$ for all $r \in (0, R]$.

\end{document}